\documentclass[11pt,leqno]{amsart}
\usepackage{amssymb, amsmath}
\usepackage{pstricks,pst-node}

\let\oldsection=\section

\makeatletter
\renewcommand{\@seccntformat}[1]{\bf\@nameuse{the#1}.\quad}
\renewcommand\section{\@startsection{section}{1}%
            \z@{.7\linespacing\@plus\linespacing}{.5\linespacing}%
            {\normalfont\bfseries \boldmath}}
\renewcommand\subsection{\@startsection{subsection}{2}%
            \z@{.5\linespacing\@plus.7\linespacing}{-.5em}%
            {\normalfont\bfseries \boldmath}}
\renewcommand\subsubsection{\@startsection{subsubsection}{3}%
            \z@{.3\linespacing\@plus.5\linespacing}{-.5em}%
            {\normalfont\bfseries \boldmath}}
\makeatother

\theoremstyle{plain}
\newtheorem*{thm}{Theorem}

\newtheorem*{lem}{Lemma}
\newtheorem*{prop}{Proposition}

\theoremstyle{definition}
\newtheorem*{defn}{Definition}
\newtheorem*{rem}{Remark}
\newtheorem*{ex}{Example}

\numberwithin{equation}{subsection}

\newcounter{listequation}

\def\note#1{{\small\tt <<#1>>}}  
\def\note#1{}              

\def\N{\mathbb N}
\def\:{\colon}
\def\a{\alpha}

\def\g{\gamma}
\def\d{\delta}
\def\e{\varepsilon}

\def\la{\lambda}

\def\cal{\mathcal}

\def\cal{\mathcal}

\def\phi{\varphi}

\def\C{\mathbb C}
\def\CC{\mathbb C}

\def\RR{\mathbb R}
\def\N{\mathbb N}

\def\ZZ{\mathbb Z}
\def\D{\mathcal D}
\def\:{\colon}

\def\SW{{}^S W}

\def\SWJ{{}^S W^J}

\def\JWw0S{{}^J W^{-w_0 S}}
\def\w0JWS{{}^{-w_0 J} W^S}
\def\fg{{\mathfrak g}}
\def\fh{{\mathfrak h}}
\def\fb{{\mathfrak b}}
\def\fp{{\mathfrak p}}

\def\fm{{\mathfrak m}}

\def\fu{{\mathfrak u}}
\def\fq{{\mathfrak q}}
\def\fl{{\mathfrak l}}

\def\fa{{\mathfrak a}}

\def\Hom{\operatorname{Hom}}
\def\Ext{\operatorname{Ext}}
\def\Tor{\operatorname{Tor}}

\def\O{\mathcal O}
\def\Pal{\text{Pal}}

\def\TJJ'{T_J^{J'}}
\def\TJ'J{T_{J'}^J}
\def\reg{\operatorname{reg}}

\def\Label{\label}

\begin{document}
\title[Kostant modules in blocks of category ${\cal O}_{S}$]
{Kostant modules in blocks of category ${\cal O}_{S}$}

\author{Brian D. Boe}
\thanks{Research of the first author partially supported by NSA grant H98230-04-1-0103}
\address{Department of Mathematics \\
                      University of Georgia\\ Athens, Georgia 30602  }
\email{brian@math.uga.edu}

\author{Markus Hunziker}
\address{Department of Mathematics \\
            Baylor University\\ Waco, Texas 76798 }
\email{Markus\_Hunziker@baylor.edu}

\date{\today}

\subjclass{}

\keywords{}

\dedicatory{}

\begin{abstract} 
In this paper the authors investigate infinite-dimensional representations $L$
in blocks of the relative (parabolic) category $\O_{S}$ for a complex simple Lie algebra, having the property that the cohomology of the nilradical with coefficients in $L$ ``looks like'' the cohomology with coefficients in a finite-dimensional module, as in Kostant's theorem. 
A complete classification of these ``Kostant modules'' in regular blocks for maximal parabolics in the simply laced types is given. A complete classification 
is also given in arbitrary (singular) blocks for Hermitian symmetric categories.
\end{abstract}

\maketitle

\parskip=2pt

\section{Introduction}
\Label{S:Intro}

\subsection{} 
Let $\fg$ be a complex simple Lie algebra with a standard parabolic subalgebra $\fp_{S}=\fm_{S}\oplus\fu_{S}$, where $\fu_{S}$ is the nilradical and $\fm_{S}$ the  Levi subalgebra, corresponding to a subset $S$ 
of the simple roots of $\fg$. Let $W$ be the Weyl group of $\fg$ and 
$W_{S}$ the parabolic subgroup of $W$ corresponding to $S$.
In \cite{Kos:61}, Kostant proved that 
for any finite-dimensional simple $\fg$-module $E$, 
the Lie algebra cohomology $H^{*}(\fu_{S}, E)$ 
is multiplicity-free as an $\fm_{S}$-module, and that the 
decomposition is described by the graded poset ${}^{S}W$ of minimal length 
coset representatives of $W_{S}\backslash W$.
Explicitly, if $E=E_{\lambda}$ is a simple finite-dimensional
$\fg$-module of highest weight $\lambda$, then as an
$\fm_{S}$-module,
\begin{equation} \Label{E:Kostant}
    H^i(\fu_{S},E_\lambda)\simeq 
    \bigoplus_{\begin{smallmatrix}y \in {}^{S}\!W \\
    l(y)=i\end{smallmatrix}}  
     F_{y(\lambda+\rho)-\rho},
\end{equation}
where $l(y)$ is the length of $y\in {}^{S}W$  and 
$F_{y(\lambda+\rho)-\rho}$ is a simple $\fm_{S}$-module of highest weight 
$y(\lambda+\rho)-\rho$. Here $\rho$ is half the sum of the positive roots 
of $\fg$ as usual.

\subsection{}
In this paper we investigate infinite-dimensional simple highest weight
modules of  $\fg$ that have $\fu_{S}$-cohomology formulas analogous
to Kostant's formula \eqref{E:Kostant}. 
We will consider these ``Kostant modules'' in the relative (parabolic) Bernstein-Gelfand-Gelfand
category $\O_{S}$ that was introduced by Rocha-Caridi \cite{Roc:80}. 
The relative category $\O_{S}$ can be decomposed 
into a direct sum of  blocks, each block being a full subcategory consisting of modules with the same generalized central character. Each block contains
only finitely many simple modules.
Regular blocks, corresponding to regular central characters, 
are the easiest to describe. The simple modules in a regular block 
are in one-to-one correspondence with the poset ${}^{S}W$ (see Section~\ref{S:Prelim}). Let $L_{w}$ denote the simple module in a regular block corresponding to $w\in {}^{S}W$. Then $L_{w}$ is a Kostant module if and only if the relative (parabolic) Kazhdan-Lusztig polynomials ${}^{S}P_{x,w}$ are all zero or one. We will use this characterization to give a complete classification of Kostant modules in regular
blocks for maximal parabolics in the simply laced types in terms of Dynkin diagrams (see Sections~\ref{S:standardKostant} and \ref{S:ADEmaxreg}).

The situation is more complicated for singular blocks. 
It is still true that the simple modules in singular blocks are in one-to-one correspondence with a certain poset, but this poset is not a ranked poset in general. This fact makes it more difficult to characterize Kostant modules
in general singular blocks (see Section~\ref{S:Def} for the  precise definition of a Kostant module in singular blocks). However, in certain nice special cases,
for example in the Hermitian symmetric cases, the poset parameterizing the simple
modules in a singular block is always isomorphic to a poset of the form 
${}^{S'}W'$ and the Kostant modules correspond to Kostant modules in a regular block of category $\O_{S'}$ for some simple Lie algebra $\fg'$ of smaller rank.
Using results of Enright and Shelton we will give a complete classification
of Kostant modules in singular blocks for Hermitian symmetric pairs
in terms of Dynkin diagrams (see Section~\ref{S:HSsing}). While studying 
singular blocks in non-Hermitian symmetric cases, we discovered a new 
ordering (i.e., different from the usual Bruhat ordering) on the posets parameterizing the simple modules (see Section~\ref{S:General}). We expect that this new ordering will be useful in the study of singular blocks of 
category $\O_{S}$ in general (not just in the study of Kostant modules).

\subsection{}
One motivation to study Kostant modules is because of their connection to rationally smooth Schubert varieties (see Section~\ref{SS:SchubertCorresp}). Let $G$ be a connected complex algebraic group with Lie algebra $\fg$ and let $B\subset P\subset G$ be the connected closed subgroups corresponding to $\fb\subset\fp\subset \fg$. The $B$-orbits in $G/P$ are parameterized by $\SW$, and the closure of the orbit parameterized by $w\in\SW$ is a Schubert variety $Y_w$. By work of Kazhdan-Lusztig and Deodhar, $Y_w$ is rationally smooth if and only if every relative Kazhdan-Lusztig polynomial ${}^S P_{x,w}$ is zero or one. This is the same condition for the simple module $L_w$ in a regular block of $\O_S$ to be a Kostant module. Our uniform classification of the Kostant modules in regular blocks for maximal parabolics in simply-laced type, and for arbitrary Hermitian symmetric pairs, via subdiagrams of the Dynkin diagram, gives the simplest and most explicit description yet obtained of the rationally smooth Schubert varieties in those cases.

\subsection{}
Another of our primary motivations to study and classify Kostant modules was to understand which simple highest weight modules admit a Bernstein-Gelfand-Gelfand (BGG) resolution. It turns out that in regular blocks of category $\O_{S}$, the simple modules that admit a BGG resolution are precisely the Kostant modules (see Section~\ref{S:BGG}). Again using results of
Enright and Shelton it also follows that every Kostant module in a singular block for a Hermitian symmetric pair admits a BGG resolution. Of particular
interest are BGG resolutions of unitarizable highest weight modules. It was already observed in \cite{EnHu:04a}, that the BGG resolutions of certain
unitarizable highest weight modules can be interpreted as minimal free resolutions of coordinate rings of determinantal varieties. The calculation
of these resolutions using the structure of Kostant modules is 
remarkably simple. To illustrate the power of this technique,
we calculate the minimal free resolution of the coordinate ring 
of the closure of the minimal $E_{6}(\CC)$-orbit in the $27$-dimensional representation (see Section~\ref{S:HSsing}).

 \subsection{\it Acknowledgements} The notion of a Kostant module is not entirely new. 
 In \cite{Col:85}, Collingwood 
 introduced the notion of Kostant module for Harish-Chandra modules
 and classified what in our terminology would be  Kostant modules 
 with regular infinitesimal character for Hermitian symmetric pairs.  
 Our result in Section~\ref{S:ADEmaxreg} is a generalization of Collingwood's
 result. Also in the Hermitian symmetric
 pair setting, Enright in \cite{Enr:88} proved a 
 $\fu$-cohomology formula analogous to Kostant's formula for all unitary  highest weight modules (with not necessarily regular infinitesimal character)
 by using methods that were developed earlier by Enright and Shelton in \cite{EnSh:87}. Our result in Section~\ref{S:HSsing} is a generalization of 
Enright's result.

 The authors thank Daniel Nakano and Jonathan Kujawa for helpful discussions.

\section{Notation and Preliminaries}
\Label{S:Prelim}

\subsection{\it Notation.}
\begin{figure}[ht]
\centering
\begin{pspicture}(-.3,-.35)(0.3,0.5)
$
\cnode*(-1,0){.07}{a}
\cnode*(0,0){.07}{b}
\cnode*(1,0){.07}{c}
{\psset{linewidth=.5pt,labelsep=8pt}
\ncline{a}{b}
\ncline{b}{c}
\psline(-0.15,0.15)(0.15,-0.15)
\psline(-0.15,-0.15)(0.15,0.15)
{\small
\uput[u](-1,0){a}
\uput[u](0,0){b}
\uput[u](1,0){c}
}
}
$
\end{pspicture}
\caption{Dynkin diagram for $(A_{3},A_{1}\times A_{1})$}
\label{F:A3x2}
\end{figure}
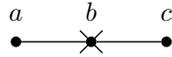

In the following let $\fg$ be a (finite-dimensional) 
complex simple Lie algebra with fixed
Cartan subalgebra $\fh$ and Borel subalgebra $\fb$ containing $\fh$.
Let $\Phi\subset \fh^{*}$ denote the root system of $(\fg,\fh)$. 
For $\alpha \in \Phi$,
let $\fg^{\alpha}$ denote the root subspace of $\fg$ corresponding to $\alpha$.
If $\fa$ is an $\fh$-invariant subspace of $\fg$, let 
$\Phi(\fa)=\{\alpha\in \Phi : \fg^{\alpha} \subset \fa\}$.
Let $\Phi^{+}=\Phi(\fb)$ be the set of positive roots and let $\Delta$ be the
subset of simple roots in $\Phi^{+}$. Let $\rho$ be half the sum of all the 
positive roots. Every subset $S\subset \Delta$ determines in the usual way
a standard parabolic subalgebra $\fp_S=\fm_S\oplus \fu_S$ 
containing $\fb$ with Levi factor $\fm_S$ and nil radical $\fu_S$ 
such that $\Phi(\fm_S)\cap \Delta = S$ and 
$\Phi(\fu_S)\cap \Delta =\Delta\smallsetminus S$. We will frequently specify 
a parabolic subalgebra $\fp_{S}$ by giving the pair $(\Phi, \Phi_{S})$, where $\Phi_{S}=\Phi(\fm_{S})$.
Define the Dynkin diagram of the pair 
$(\fg,\fp_{S})$ as the Dynkin diagram of $\fg$ with the nodes corresponding to the simple roots in $\Delta\smallsetminus S$ crossed. By the \emph{type} of such a diagram 
we will mean the pair of types of $(\Phi,\Phi_{S})$.
For example, Figure~\ref{F:A3x2}
shows the diagram of type $(A_{3},A_{1}\times A_{1})$ of the simple Lie algebra of type $A_{3}$ and its standard parabolic subalgebra corresponding to $S=\{a,c\}$. 
 
\subsection{\it Category $\O_{S}$.}
Fix $S\subset \Delta$ and let $\fp=\fm\oplus \fu $ be the standard parabolic corresponding $S$. (We will usually omit the subscript $S$ when $S$ has been fixed.)
Let $\O_{S}$ be the category of all finitely generated $\fg$-modules
on which $\fm$ acts finitely semisimply and $\fu$ acts locally nilpotently.
If $S=\varnothing$  (i.e., if $\fp=\fb$), then
the category $\O_{S}$ is the ordinary Bernstein-Gelfand-Gelfand 
category $\O$. The most basic objects in $\O_{S}$ are the generalized 
Verma modules. Let $F_{\la}$ be a finite-dimensional 
simple $\fm$-module with highest weight $\la$. We may view
$F_{\lambda}$ as a $\fp$-module by letting the nilradical $\fu$ act trivially.
Then the generalized Verma module (GVM) with highest weight $\lambda$ is the
induced module
$$
   N_{\lambda} =U(\fg)\otimes_{U(\fp)} F_{\lambda}.
$$
The module $N_{\lambda}$ is a quotient of the ordinary Verma module 
$M_{\lambda}$ with highest weight $\lambda$. Let $L_{\lambda}$ denote 
the simple quotient of both $N_{\lambda}$ and $M_{\lambda}$.
The Verma module $M_{\lambda}$ and all of its subquotients, 
in particular $N_{\lambda}$ and $L_{\lambda}$, admit a central character 
which we denote $\chi_{\lambda}$. Recall that this means that 
$\chi_{\lambda}:Z(U(\fg))\to \CC$ is a character of the center
of $U(\fg)$ such that $z.m= \chi_{\lambda}(z)m$ for all 
$z\in Z(U(\fg))$ and all $m$ in the module.
For any character $\chi :Z(U(\fg))\to \CC$ let  
$\O_{S}^{\,\chi}$ denote the full subcategory of $\O_{S}$ consisting 
of modules $V$ such that $z-\chi(z)$ acts locally nilpotently on $V$ 
for all $z\in Z$. 
Then the category $\O_{S}$ decomposes as 
$$
\O_{S}=\bigoplus_{\chi\in Z(U(\fg))^{\wedge}}\O_{S}^{\,\chi}.
$$ 
Here $\O_{S}^{\,\chi}$ is non-zero if and only if $\chi$ is of the form
$\chi=\chi_{\lambda}$ for some $S$-dominant $S$-integral $\lambda\in \fh^{*}$.
In the following, if $\chi=\chi_{\lambda}$ we will write $\O_{S}^{\,\lambda}$
instead of $\O_{S}^{\,\chi}$. The categories $\O_{S}^{\lambda}$ are 
called the blocks of $\O_{S}$.
Let $W$ denote the Weyl group of $\fg$. As usual we write $w\cdot\lambda$
for $w(\lambda+\rho)-\rho$.
By a theorem of Harish-Chandra, $\chi_{\lambda}=\chi_{\mu}$,
and hence $\O_S^\lambda=\O_S^\mu$, if and only if $\lambda\in W\cdot \mu$.
Thus, without loss of generality, we may consider only
blocks of the form $\O_{S}^{\,\mu}$, where $\mu+\rho$ is 
anti-dominant. Furthermore, for simplicity, we will only consider 
integral weights. So, from now on we will assume that  
$\mu+\rho$ is an anti-dominant integral weight,
i.e., $(\mu+\rho,\alpha^\vee)$ is a non-positive integer for 
all $\alpha\in \Delta$. Here 
$(\underline{\ \ },\underline{\ \ })$ denotes the bilinear form  on $\fh^{*}$
(induced from the Killing form) and 
$\alpha^{\vee}=2\alpha/(\alpha,\alpha)$.

\subsection{\it Posets and parameterization of simple modules in $\O_{S}^{\mu}$}\Label{SS:posets}
Let $W_{S}$ be the subgroup of $W$ generated by the simple reflections
$s_{\a}$ with $\a\in S$. Define 
$$
{}^{S} W = \{ w \in W \mid l(s_{\a}w) = l(w) + 1\  \mbox{for all}\ \a \in S\},
$$
where $l$ is the usual length function on $W$.
Then ${}^{S}W$ is the set of minimal length coset representatives
for $W_{S}\backslash W$. There is a natural ordering on ${}^{S}W$
which is induced by the Bruhat ordering on $W$ (cf.\ \cite{Hum:90}).
We view ${}^{S}W$ as a poset (partially ordered set) 
with respect to this ordering and we will later often draw 
a Hasse diagram of ${}^{S}W$ with the least element $e$ at the bottom  
(e.g., Figures~\ref{F:D4A2}, \ref{F:F4C3}). We sometimes label an edge 
$x$---$w$ by a simple root $\a$ if $w=xs_\a$ (as in Figure \ref{F:F4C3}).

Now let $\mu\in \fh^{*}$ be such that $\mu+\rho$ is anti-dominant integral
and let $w_{S}$ be the longest element of $W_{S}$. 
Then every highest weight $\lambda$ of a highest weight module in 
$\O_{S}^{\,\mu}$ can be written in the form 
$\lambda=w_Sw\cdot\mu$ for some $w\in {}^{S}W$.
If $\mu+\rho$ is not regular then the element $w\in {}^{S}W$ is not unique.
Define the set of singular simple roots by
$$
J = \{ \alpha \in \Delta \mid (\mu+\rho ,\alpha^\vee)=0\}.
$$
Then $W_J=W_\mu=\{w\in W \mid w\cdot \mu =\mu\}$.
Define
$$
 {}^SW^J =\{ w\in {}^SW \mid  w < ws_\alpha \in  {}^SW \mbox{\ for all $\alpha\in J$}\}.
$$
In the Hasse diagram of ${}^SW$, a node that corresponds
to an element of $\SWJ$ is a node having an edge with label $\alpha$
going up from it for every $\alpha \in J$.
For now we will view ${}^SW^J$ as a poset with ordering induced from
the Bruhat ordering on ${}^{S}W$. 
(A different ordering will be introduced later in Section~\ref{S:General}.)

Every highest weight $\lambda$ of a highest weight module in 
$\O_{S}^{\,\mu}$ can be written {\it uniquely} in the form 
$\lambda=w_Sw\cdot\mu$ for some $w\in {}^{S}W^{J}$.
In the following we will write $L_{w}$ (resp.\ $N_{w}$) for 
the simple module (resp.\ GVM) in $\O_{S}^{\, \mu}$ of highest weight 
$w_Sw\cdot\mu$. We also write $F_w$ for the simple finite-dimensional $\fm$-module of highest weight $w_Sw\cdot\mu$.

\section{Definition of Kostant modules}
\Label{S:Def}

\subsection{}
In this section we define Kostant modules in regular and singular blocks of 
category $\O_{S}$. In the regular setting, we also give several equivalent characterizations of Kostant modules, which will use in our classification arguments. In regular blocks, it is quite clear what the ``correct'' 
definition should be. In singular blocks, it is however not so clear,
and the definition that we give below should be thought of as a ``working definition'' that may have to be modified in future work. 
The difficulty results from the structure of the posets $\SWJ$ that can be 
quite complicated when $J\not=\varnothing$. We will come back to this issue
in Section~\ref{S:General}.

\subsection{\it Graded posets}
Recall that a finite poset is called \emph{graded} if it is an interval
(i.e., it has a least and a greatest element) and if all maximal chains between any two elements $x<y$
have the same length. On a graded interval $[v,w]$ we then have a well-defined {\it rank function} given by $r(x)$=length of any maximal chain from $v$ to $x$.
We will write $[v,w]_{j}$ for the set of all $x$ in $[v,w]$ for which $r(x)=j$.

It is well-known (e.g., see Deodhar \cite{Deo:77}), 
that the posets of the form $\SW$ are graded and that the rank function on  
$\SW$ is the restriction of the length function $l$ on $W$.
However, the posets of the form $\SWJ$ are not graded in general.

\subsection{}
Let $\O_{S}^{\,\mu}$ be a regular or singular block of category $\O_{S}$.

\begin{defn}\label{D:Kostant}
For $w\in \SWJ$, 
we say that $L_{w}$ is a \emph{Kostant module} in $\O_{S}^{\,\mu}$
if there exists a graded interval  $[v,w]$ of  $\SWJ$ such that
as an $\fm$-module, 
$$
  H^{i}(\fu,L_{w}) \simeq \bigoplus_{x\in [v,w]_{r(w)-i}} F_{x},
$$
where $r$ is the rank function on $[v,w]$.
\end{defn}

\begin{rem}
Every finite-dimensional simple module is a Kostant module
by Kostant's theorem. Note that in our notation, a simple module 
$L_{w}$ in $\O_{S}^{\,\mu}$ is 
finite-dimensional if and only if $\mu+\rho$ is regular and
$w$ is the (unique) maximal element of $\SWJ={}^{S}W$. The interval in this
case is $[v,w]=[e,w]={}^{S}W$, and the rank function is the length function.
\end{rem}

\subsection{\it Kostant modules and Kazdan-Lusztig-Vogan polynomials} \Label{SS:KLV}
For regular blocks the definition of Kostant modules has a number of equivalent formulations. The first of these is in terms of Kazhdan-Lusztig-Vogan polynomials.
Let $\O_{S}^{\, \mu}=\O_{S}^{\reg}$ be any regular block; that is, 
$\mu+\rho$ is regular. 
For $x,w \in {}^{S}W$ define the relative 
{\it Kazdan-Lusztig-Vogan polynomial\/}  
$$
{}^{S}P_{x,w}= \sum_{i\geq 0} q^{\frac{l(w)-l(x)-i}{2}} 
\dim \Ext^{i}_{\O_{S}}(N_x, L_w).
$$
It follows from the Kazhdan-Lusztig conjectures (which are theorems in this setting; cf.~\cite{CaCo:87}) that the polynomial
${}^{S}P_{x,w}$ is equal to the parabolic Kazhdan-Lusztig polynomial
that was defined by Deodhar \cite{Deo:87}.
It is well known (and follows, for example, by a straightforward extension of the proof of \cite[Lemma~5.13]{Sch:81}, the corresponding statement in ordinary category $\O$), that
$$
    \Hom_\fm(F_x, H^i(\fu, L_w)) \simeq \Ext^{i}_{\O_{S}}(N_x,L_w)\quad 
    \mbox{as vector spaces}.
$$
It then follows immediately from the definition of a 
Kostant module in Definition~\ref{D:Kostant} (and the fact that ${}^{S}P_{e,w}\not=0$)
that the simple module $L_w$ 
in $\O_S^{\reg}$ is a Kostant module  if and only if 
$${}^{S}P_{x,w}=1\quad \text{ for every } x\in [e,w].$$

\subsection{\it Schubert varieties.}\Label{SS:SchubertCorresp}
Further characterizations of Kostant modules in regular blocks involve connections with the geometry of Schubert varieties. Let $G$ be a connected complex algebraic group with Lie algebra $\fg$ and
let $T\subset B\subset P \subset G$ be the connected closed subgroups
corresponding to  $\fh\subset\fb\subset\fp\subset \fg$.
The maximal torus $T$ acts with finitely many fixed points on 
the generalized flag variety $Y=G/P$ and 
these fixed points are naturally parameterized by the parabolic poset 
${}^{S}W$. The $B$-orbits in $Y$ are the orbits of the $T$-fixed points.
The closure of the orbit through the $T$-fixed point corresponding
to $w\in {}^{S}W$ is denoted $Y_{w}$ and is called a generalized
{\it Schubert variety}. It follows from the work of Kazhdan-Lusztig 
\cite{KaLu:79} and Deodhar \cite{Deo:87} that $Y_{w}$ is 
rationally smooth if and only if ${}^{S}P_{x,w}=1$ for every 
$x\in [e,w]$.

Thus, by the characterization in Section~\ref{SS:KLV}, we have a canonical 1-1 correspondence
$$
\{\mbox{Kostant modules in $\O_S^{\reg}$}\} 
\leftrightarrow
\{\mbox{rationally smooth Schubert varieties in $Y$}\},
$$
where $L_{w}$ corresponds to $Y_{w}$. To classify the Kostant modules in regular blocks we can then use any of the many known tests for rational smoothness of 
Schubert varieties (see Billey-Lakshmibai \cite{BiLa:00} for an overview). 
A test that will be useful later is the following, due to Carrell-Peterson \cite{Car:94}. For $w\in {}^{S}W$, define the Poincar\'e
polynomial $P_{w}(t)$ by
$$
  P_{w}(t)=\sum_{v\leq w} t^{l(v)}.
$$
Then $$Y_{w} \text{ is rationally smooth if and only if } P_{w}(t) \text{ is palindromic.}$$ We will frequently use this criterion for a simple module $L_w$ in $\O_S^{\reg}$ to be a Kostant module.

\begin{ex} Consider $(\Phi,\Phi_{S})$ of type $(D_{4},A_{2})$.
Figure~\ref{F:D4A2} shows the Hasse diagram of the poset ${}^{S}W$. The circled nodes of the Hasse diagram  correspond to Kostant modules, or equivalently, smooth Schubert varieties. (In simply-laced types, a result of Peterson (see  \cite{CaKu:03}) says that a Schubert variety is rationally smooth if and only if it is smooth.)
The reader may want to verify that the circled nodes are precisely the
nodes for which the corresponding Poincar\'e polynomial is palindromic. 

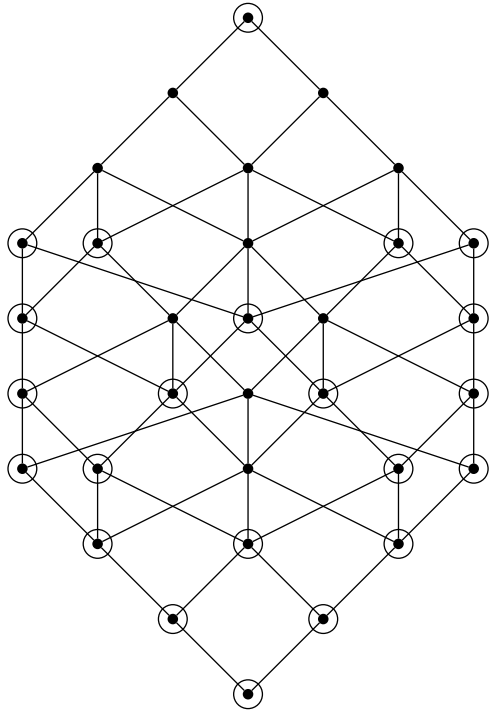
\begin{figure}[ht]
\centering
\begin{pspicture}(-.3,-.35)(.2,10)
\psset{linewidth=.5pt,labelsep=8pt,nodesep=0pt}
\small
$
\cnode*(0,0){.07}{a1}
\pscircle(0,0){.2}
\cnode*(-1,1){.07}{a2}
\pscircle(-1,1){.2}
\cnode*(1,1){.07}{a3}
\pscircle(1,1){.2}
\cnode*(-2,2){.07}{a4}
\pscircle(-2,2){.2}
\cnode*(0,2){.07}{a5}
\pscircle(0,2){.2}
\cnode*(2,2){.07}{a6}
\pscircle(2,2){.2}
\cnode*(-3,3){.07}{a7}
\pscircle(-3,3){.2}
\cnode*(-2,3){.07}{a8}
\pscircle(-2,3){.2}
\cnode*(0,3){.07}{a9}
\cnode*(2,3){.07}{a11}
\pscircle(2,3){.2}
\cnode*(3,3){.07}{a10}
\pscircle(3,3){.2}
\cnode*(-3,4){.07}{a12}
\pscircle(-3,4){.2}
\cnode*(-1,4){.07}{a13}%
\pscircle(-1,4){.2}
\cnode*(0,4){.07}{a14}
\cnode*(1,4){.07}{a15}%
\pscircle(1,4){.2}
\cnode*(3,4){.07}{a16}
\pscircle(3,4){.2}
\cnode*(-3,5){.07}{a17}
\pscircle(-3,5){.2}
\cnode*(-1,5){.07}{a18}
\cnode*(0,5){.07}{a19}
\pscircle(0,5){.2}
\cnode*(1,5){.07}{a20}
\cnode*(3,5){.07}{a21}
\pscircle(3,5){.2}
\cnode*(-3,6){.07}{a23}
\pscircle(-3,6){.2}
\cnode*(-2,6){.07}{a22}
\pscircle(-2,6){.2}
\cnode*(0,6){.07}{a24}
\cnode*(2,6){.07}{a25}
\pscircle(2,6){.2}
\cnode*(3,6){.07}{a26}
\pscircle(3,6){.2}
\cnode*(-2,7){.07}{a27}
\cnode*(0 ,7){.07}{a28}
\cnode*(2,7){.07}{a29}
\cnode*(-1 ,8){.07}{a30}
\cnode*(1 , 8){.07}{a31}
\cnode*(0 ,9){.07}{a32}
\pscircle(0,9){.2}
\ncline{a1}{a2}
\ncline{a1}{a3}
\ncline{a2}{a4}
\ncline{a2}{a5}
\ncline{a3}{a5}
\ncline{a3}{a6}
\ncline{a4}{a7}
\ncline{a4}{a8}
\ncline{a4}{a9}
\ncline{a5}{a8}
\ncline{a5}{a9}
\ncline{a5}{a11}
\ncline{a6}{a9}
\ncline{a6}{a11}
\ncline{a6}{a10}
\ncline{a7}{a12}
\ncline{a7}{a14}
\ncline{a8}{a12}
\ncline{a8}{a13}
\ncline{a9}{a13}
\ncline{a9}{a14}
\ncline{a9}{a15}
\ncline{a11}{a15}
\ncline{a11}{a16}
\ncline{a10}{a14}
\ncline{a10}{a16}
\ncline{a12}{a17}
\ncline{a12}{a18}%
\ncline{a13}{a17}%
\ncline{a13}{a18}
\ncline{a13}{a19}
\ncline{a14}{a18}
\ncline{a14}{a20}
\ncline{a15}{a19}
\ncline{a15}{a20}
\ncline{a15}{a21}%
\ncline{a16}{a20}%
\ncline{a16}{a21}
\ncline{a17}{a23}
\ncline{a17}{a22}
\ncline{a18}{a22}
\ncline{a18}{a24}
\ncline{a19}{a23}%
\ncline{a19}{a24}
\ncline{a19}{a26}%
\ncline{a20}{a24}
\ncline{a20}{a25}
\ncline{a21}{a25}
\ncline{a21}{a26}
\ncline{a23}{a27}
\ncline{a22}{a27}
\ncline{a22}{a28}%
\ncline{a24}{a27}%
\ncline{a24}{a28}
\ncline{a24}{a29}%
\ncline{a25}{a28}%
\ncline{a25}{a29}
\ncline{a26}{a29}
\ncline{a27}{a30}
\ncline{a28}{a30}
\ncline{a28}{a31}
\ncline{a29}{a31}
\ncline{a30}{a32}
\ncline{a31}{a32}
$
\end{pspicture}
\caption{Kostant modules in a regular block for 
$(D_{4},A_{2})$}
\label{F:D4A2}
\end{figure}

\end{ex}

\section{Generalized Bernstein-Gelfand-Gelfand resolutions}
\Label{S:BGG}

\subsection{}
In this section we prove that every Kostant module in a regular block 
$\O_{S}^{\reg}$ admits a resolution in terms of sums of generalized Verma modules in $\O_{S}^{\reg}$. The technical part of the proof was already 
given in \cite{EnHu:04b}. 

\subsection{\it Definition of the BGG complex}
The construction of the resolution of a Kostant module $L_{w}$ 
in $\O_{S}^{\reg}$ is combinatorial and completely analogous to 
the construction of the resolution of a finite dimensional simple module 
that was given by Lepowsky \cite{Lep:77} 
(and originally by Bernstein-Gelfand-Gelfand \cite{BGG:75}
for $S=\varnothing$). For the convenience of the reader, we repeat some
of the details. By \cite[Prop.~3.7]{Lep:77}, for every arrow  $x\rightarrow y$ in ${}^{S}W$ there exists a nonzero $\fg$-module map $f_{x,y}: N_{x}\rightarrow N_{y}$, which lifts to the standard map between the (ordinary) Verma modules
having the same highest weights as $N_{x}$ and $N_{y}$, respectively.
Recall that a quadruple $(w_{1},w_{2},w_{3},w_{4})$ of elements
in $W$ is called a square if $w_{2}\not=w_{3}$ and 
$w_{1}\rightarrow w_{2} \rightarrow w_{4}$
and $w_{1}\rightarrow w_{3} \rightarrow w_{4}$. 
By \cite[Lemma~10.4]{BGG:75}, 
it is possible to assign to each arrow $x\rightarrow y$ in $W$
a number $\varepsilon_{x,y}=\pm 1$ such that for every square,
the product of the four numbers assigned to the sides of the square
is $-1$. Now for any $x,y\in {}^{S}W$ with  $l(y)=l(x)+1$,  define
a $\fg$-module map $h_{x,y}:N_{x} \rightarrow N_{y}$ by
$$
h_{x,y}=
\begin{cases}
\varepsilon_{x,y}f_{x,y}  & \mbox{if $x<y$};\\
0 & \mbox{otherwise}.
\end{cases}
$$
Let $x,z \in {}^{S}W$ with $l(z)=l(x)+2$. 
Then the number of elements $y\in {}^{S}W$ such that 
$x\rightarrow y\rightarrow z$ is either 
zero, one, or two. If there is one, then $f_{y,z}\circ f_{x,y}=0$ 
and hence $h_{y,z}\circ h_{x,y}=0$. 
If there are two, i.e., if $x\rightarrow y\rightarrow z$ and $x\rightarrow y'\rightarrow z$ with $y\not=y'$, then 
$h_{z,y}\circ h_{x,y} = - h_{z,y'}\circ h_{x,y'}$.
These observations due to Lepowsky can now be used to construct  
a complex for {\it any\/} simple module $L_{w}$ in $\O_{S}^{\reg}$.
For  $0\leq i \leq l(w)$, define 
$$
  C_{i} =\bigoplus_{x\in [e,w]_{l(w)-i}} N_{x}
$$
and for $1\leq i \leq l(w)$, define $d_{i} : C_{i} \rightarrow C_{i-1}$
as the matrix of maps $d_{i}=(h_{x,y})$, where $x\in [e,w]_{l(w)-i}$ and 
$y\in [e,w]_{l(w)-(i-1)}$. Furthermore, let 
$d_{0}:C_{0}=N_{w}\rightarrow L_{w}$ be the canonical quotient map.

\begin{lem}
The sequence
$$   0\rightarrow C_{l(w)} \rightarrow  \cdots \rightarrow
  C_{1} \rightarrow C_{0}
  \rightarrow L_{w} \rightarrow 0
$$
is a complex, i.e., $d_{i-1}\circ d_{i}=0$ for $1\leq i\leq l(w)$.
Furthermore, for every $x\in [e,w]_{l(w)-i}$, the restriction of 
$d_{i}$ to $N_{x}$ is nonzero.
\end{lem}

\begin{proof}
Note that if $x,z \in [e,w]$ and $y\in {}^{S}W$ with 
$x\rightarrow y \rightarrow z$, then $y\in [e,w]$. The observations above then 
immediately imply that $d_{i-1}\circ d_{i}=0$ for $2 \leq i \leq l(w)$.
Thus it remains to show that $d_{0} \circ d_{1}=0$. This follows since for 
every $x\rightarrow w$ in ${}^{S}W$, the image of 
$f_{x,y}: N_{x} \rightarrow N_{w}$ is contained in the radical of 
$N_{w}$, which is equal to the kernel of the quotient map $N_{w}\rightarrow L_{w}$.
\end{proof}

\subsection{\it Existence of BGG resolutions} We now have the following result.

\begin{thm} 
Let $L_{w}$ be a simple module in $\O_{S}^{\reg}$. Then
the complex $0\rightarrow C_{l(w)} \rightarrow  \cdots \rightarrow
  C_{1} \rightarrow C_{0}
  \rightarrow L_{w} \rightarrow 0$
is exact if and only if $L_{w}$ is a Kostant module. 
\end{thm}

\begin{proof}
Suppose that $L_{w}$ is a Kostant module. By the lemma above, $L_{w}$ 
is then a generalized Kostant module in the sense of \cite[2.7]{EnHu:04b}, and
\cite[Thm.~2.8]{EnHu:04b} shows that the complex $0\rightarrow C_{l(w)} \rightarrow  \cdots \rightarrow
  C_{1} \rightarrow C_{0}
  \rightarrow L_{w} \rightarrow 0$
is exact.

To prove the converse suppose that $L_{w}$ is a simple module such that the 
complex $0\rightarrow C_{l(w)} \rightarrow  \cdots \rightarrow
  C_{1} \rightarrow C_{0}
  \rightarrow L_{w} \rightarrow 0$
is exact.
Let $\fg=\fu^{-}\oplus\fm\oplus\fu$ be the triangular decomposition 
of $\fg$ given by $S$. It follows from the PBW Theorem that as an $\fm$-module,
the generalized Verma module $N_{x}$ is isomorphic to $N_{x}\simeq U(\fu^{-})\otimes_{\CC} F_{x}\simeq S(\fu^{-})\otimes_{\CC} F_{x}$ for all $x\in {}^{S}W$.  
The Killing form induces an isomorphism between $S(\fu)$ and $S(\fu^{-})$. 
This in turn induces an $\fm$-equivariant isomorphism between $H^{i}(\fu,L_{w})=\Ext^{i}_{\fu}(\CC,L_{w})$ and $H_{i}(\fu^{-},L_{w})=\Tor_{i}^{\fu^{-}}(\CC,L_{w})$. (We have a canonical isomorphism $H^{i}(\fu,L_{w})\simeq H_{i}(\fu ,L_{w}^{*}{)}^{*}$, where $L_{w}^{*}$ denotes the dual of the $\fg$-module $L_{w}$.
Via the Killing form, we have $H_{i}(\fu ,L_{w}^{*}{)}^{*}\simeq 
H_{i}(\fu^{-},L_{w})$.) Recall that $\Tor_{i}^{\fu^{-}}(\CC,L_{w})$
can be computed by any projective $U(\fu^{-})$-module resolution of $L_{w}$. 
By the remark above, the resolution $0\rightarrow C_{l(w)} \rightarrow  \cdots \rightarrow
  C_{1} \rightarrow C_{0}
  \rightarrow L_{w} \rightarrow 0$
is a free $U(\fu^{-})$-resolution of $L_{w}$. Furthermore, this resolution
is $\fm$-equivariant. By applying the functor
$\CC\otimes_{U(\fu^{-})}\underline{\quad}\ $ to the complex 
$0\rightarrow C_{l(w)} \rightarrow  \cdots \rightarrow
  C_{1} \rightarrow C_{0}\rightarrow 0$
and by taking the $i$-th homology of the resulting complex
we then obtain that as an $\fm$-module,
$$
H_i(\fu^-,L_{w})\simeq \bigoplus_{x\in [e,w]_{l(w)-i}} F_{x}.
$$
Since $H_i(\fu^-,L_{w})\simeq H^{i}(\fu,L_{w})$, it follows that 
$L_{w}$ is a Kostant module.
\end{proof}

\section{Standard Kostant Modules in Regular Blocks}
\Label{S:standardKostant}

\subsection{}
In Section~\ref{S:ADEmaxreg}, we will classify the Kostant modules in regular
blocks for any maximal parabolic in simply laced type in terms of subdiagrams
of the Dynkin diagram of the pair. In this section, we explain, for a general
parabolic, how to associate Kostant modules to subdiagrams. The
Kostant modules that arise in this way we will call {\it standard\/} Kostant
modules.  In Section~\ref{S:ADEmaxreg} we will then show that for maximal
parabolics in simply laced type, every Kostant module in a regular  block
is a standard Kostant module.

\subsection{\it Standard Kostant modules}\Label{SS:standardKostant}
Fix a standard parabolic subalgebra $\fp=\fp_{S}\subset \fg$ and let 
$P\subset G$ be the parabolic subgroup corresponding to $\fp$. 
To any standard parabolic subalgebra $\fq=\fp_{I}\subset\fg$, where $I\subset\Delta$, we can
associate a smooth Schubert variety in $G/P$ as follows. 
Let $Q\subset G$ be the parabolic subgroup corresponding to $\fq$ and 
let $Q=LU$ be the Levi decomposition of $Q$. Then $L\cap B$ is a Borel 
subgroup of the reductive group $L$ and $L\cap P$ is a 
parabolic subgroup of $L$. Furthermore, $B=(L\cap B)U$.
Now consider the inclusion $L/(L\cap P)\hookrightarrow G/P$, which in fact
is a closed embedding of smooth projective varieties.

\begin{lem} 
Under the closed embedding  $L/(L\cap P)\hookrightarrow G/P$,
the image of any  $(L\cap B)$-orbit in $L/(L\cap P)$ is a 
$B$-orbit in $G/P$. In particular, the image of $L/(L\cap P)$ is 
a smooth Schubert variety in $G/P$.
\end{lem}

\begin{proof}
Let $x\in L$ and $b\in B$. Since $B=(L\cap B)U$, we may write $b=b_{L}u$
with $b_{L}\in L\cap B$ and $u\in U$.  
Then, since $U$ is normalized by $L$, 
we have $bx=b_{L}ux=b_{L}xu'$ for some $u'\in U$.
Since $U\subset P$ (in fact, $U\subset B\subset P$), 
it follows that the $B$-orbit of $xP$ is the image
of the $(L\cap B)$-orbit of $x(L\cap P)$.
\end{proof}

\begin{rem}
The claim of the lemma is true if we replace
$L$ by the commutator  group $L'=[L,L]$, which is 
connected  semisimple. As before, $L'\cap B$ is a Borel 
subgroup of $L'$ and $L'\cap P$ is a parabolic subgroup of $L'$.
Furthermore, since $L=Z(L)L'$, where $Z(L)$ is the center of $L$, it follows
that the inclusion $L'\subset L$ induces an isomorphism 
$L'/(L'\cap P)=L/(L\cap P)$.\end{rem}

Let $Y_{w}$ be the Schubert variety that is the image of $L/(L\cap P)$
under the embedding  $L/(L\cap P)\hookrightarrow G/P$.
In light of Section~\ref{SS:SchubertCorresp}, there is then a Kostant module 
$L_w$ associated to $Y_{w}$. 
Note that $L_w$ is in fact a Kostant module because 
$[e,w]\simeq {}^{S'}W'={}^{S\cap I}W_{I}$, where $W'=W_{I}$ is the 
Weyl group of $L$ and $S'=S\cap I\subset I$ is the subset of the simple roots 
of $L$ that defines the parabolic subgroup $L\cap P \subset L$. 
The element $w$ is the maximal element of ${}^{S\cap I}W_{I}$, viewed as
an element of ${}^{S}W$.  In the following, to emphasize the dependence on $I$, we will write $w=\phi(I)$. The Kostant modules in $\O_{S}^{\reg}$
of the form $L_{\phi(I)}$ we will call {\it standard\/} Kostant modules.


\subsection{\it Subdiagrams}\Label{SS:subdiagrams}
Let $\D$ be the Dynkin diagram of the pair $(\fg,\fp)$. 
The Dynkin diagram of the pair $(\fl, \fl\cap \fp)$,
where $\fl$ is the Levi subalgebra of $\fq=\fp_{I}$, 
may be viewed as a subdiagram of $\D$. More precisely, let $\D(I)$ denote the subdiagram of $\D$ with nodes corresponding to the simple roots in $I$ and the crossed nodes corresponding to the simple roots in $I\smallsetminus S$. Then $\D(I)$ is the diagram of the pair  $(\fl, \fl\cap \fp)$.
Clearly, any subdiagram of $\D$ is of the form $\D(I)$. Furthermore, the subdiagram $\D(I)$ determines the set $I$. Thus we get a map from the set of subdiagrams of $\D$ to the set of Kostant modules in $\O_{S}^{\reg}$ given by 
$\D(I) \mapsto L_{\phi(I)}$. In the following we will distinguish a 
certain subset of the set of subdiagrams of $\D$, such that the map 
restricted to this subset is one-to-one.

Suppose that $I\not=\varnothing$ and let $I=I_{1}\cup \cdots \cup I_{r}$ be the partition of $I$ corresponding to the decomposition $\D(I)=\D(I_{1})\cup \cdots \cup \D(I_{r})$ of the diagram $\D(I)$ into connected components.
Then $\phi(I)=\phi(I_{1})\cdots \phi(I_{r})$.
We say that $\D(I_{j})$ is an {\it $S$-trivial component} of $\D(I)$
if $\phi(I_{j})=e$. Clearly, $\D(I_{j})$ is an $S$-trivial component of $\D(I)$
if and only if $I_{j}\subset S$, i.e., if and only if $\D(I_{j})$
does not contain any crossed nodes. If $I=\varnothing$ we let $\D(I)$ be the empty subdiagram of $\D$, which (vacuuously) has no $S$-trivial components. With this terminology in 
place, we have the following result.

\begin{prop}\Label{P:subdiagtoKostant}
We have an injection 
$$
\{\mbox{subdiagrams of $\D$ with no $S$-trivial components}\} \hookrightarrow
\{\mbox{Kostant modules in $\O_{S}^{\reg}$}\}
$$
given by $\D(I) \mapsto L_{\phi(I)}$.
\end{prop}

\begin{rem}
If $S^{c}=\{\alpha\}$, i.e., if $\fp$ is a maximal parabolic subalgebra of $\fg$, then the non-empty subdiagrams of $\D$ with no $S$-trivial components are the connected subdiagrams containing the node $\alpha$. For simplicity of stating results in this setting, we shall regard the empty subdiagram of $\D$ as belonging to the set of connected subdiagrams containing $\a$.
\end{rem}

The key in the proof of the proposition is the following lemma. We
assume that the lemma is well known to experts, but since we couldn't find a proof in the literature, we provide one here.

\begin{lem}\Label{L:redexpression}
Suppose that $S \not=\Delta$ and let ${}^{S}w$ be the maximal element
of ${}^{S}W$. Then any reduced expression for ${}^{S}w$ 
involves all simple reflections. 
\end{lem}

\begin{proof}
Fix $\a_p\in\Delta\smallsetminus S$. For the purposes of this proof only, number the simple roots $\a_1,\dots,\a_p,\dots,\a_m,\a_{m+1},\dots,\a_n$ where $\a_1,\dots,\a_m$ is a maximal chain containing $\a_p$ in the Dynkin diagram, and $\a_{m+1},\dots,\a_n$ is the other branch (if any) of the diagram, so that $\a_1, \a_m$, and $\a_n$ (if $m<n$) are the end nodes of the diagram. Consider the Coxeter element
$$
w:=s_p s_{p-1} \dots s_1 s_{p+1}s_{p+2}\dots s_m s_{m+1} \dots s_n.
$$
We claim that: the given expression for $w$ is reduced; every reduced expression for $w$ begins with $s_p$ and involves every simple reflection exactly once; and $w\in\SW$. 

First, any product of distinct simple reflections is necessarily reduced (cf.\ \cite[Exercise 1.13]{Hum:90}). If not, one could find an initial expression $s_{i_1} \dots s_{i_{k-1}} s_{i_k}=x s_{i_k}$ with $l(x s_{i_k}) < l(x)$; i.e., $x(\a_{i_k}) < 0$. On the other hand, by the formula for a reflection, the coefficient of $\a_{i_k}$ in $x(\a_{i_k})$, when written as a linear combination of simple roots, is 1 (since $\a_{i_1},\dots,\a_{i_{k-1}}$ are distinct from $\a_{i_k}$). This contradiction proves the first claim.

Since the set of simple reflections appearing in a fixed reduced expression for any element $x\in W$ depends only on $x$ and not on the reduced expression \cite[Cor.~5.10(b)]{Hum:90}, every reduced expression for $w$ involves every simple reflection (exactly once, by length considerations).

Any other reduced expression for a Coxeter element $w$ can be obtained from the given one by applying the ``commuting relations'' $s_i s_j = s_j s_i$ \cite[Chap.\ 4, Sec.\ 1, Prop.\ 5 and Exercise 13]{Bou:68}. If $j\ne p$, then there is an $s_i$ occuring to the left of $s_j$ in the given expression for $w$ with $s_i s_j \ne s_j s_i$. This means that $s_i$ will occur to the left of $s_j$ in every reduced expression for $w$, and in particular, no reduced expression can begin with $s_j$. By the characterization of $\SW$ as those elements of $W$ all of whose reduced expressions begin with a reflection corresponding to a simple root not in $S$, it follows that $w\in\SW$. This proves the claim.

Since ${}^Sw$ is the greatest element in the Bruhat order of $\SW$, we have $w\le {}^Sw$. By Deodhar's subexpression condition \cite[Thm.~5.10]{Hum:90}, some reduced expression for $w$ occurs as a subexpression of any given reduced expression for ${}^Sw$. By the claim, this proves the lemma.
\end{proof}

\noindent{\it Proof of Proposition.}
Let $\D(I)=\D(I_{1})\cup \cdots \cup \D(I_{r})$ be the decomposition of the (non-empty) subdiagram $\D(I)$ of $\D$ into connected components and assume
that none of the $\D(I_{j})$ is $S$-trivial. Then, by definition, $S\cap I_{j}\not=I_{j}$ for all $j$. Recall that $\phi(I_{j})$ is the maximal element of ${}^{S\cap I_{j}}W_{I_{j}}$. By the lemma above,
any reduced expression for $\phi(I_{j})$ involves all simple reflections $s_{\a}$, $\a\in I_{j}$.
Since $I=I_{1}\cup \cdots \cup I_{r}$ and $I_{j}\perp I_{j'}$ for $j\not=j'$,
this implies that any reduced expression for $\phi(I)=\phi(I_{1})\cdots\phi(I_{r})$ involves all simple reflections $s_{\a}$, $\a\in I$. Thus, $I$ is uniquely 
determined by $\phi(I)$ and the proposition follows.
\qed

\begin{ex} Consider $(\Phi, \Phi_{S})=(F_{4},C_{3})$. Figure~\ref{F:F4C3}
shows the Hasse diagram of ${}^{S}W$ with the nodes corresponding to Kostant modules circled as before. There are a total of eight Kostant modules, but only five of these eight are standard Kostant modules. The subdiagrams $\D(I)$ corresponding to $I=\varnothing$,
$\{a\}$, $\{a,b\}$, $\{a,b,c\}$, and $\Delta=\{a,b,c,d\}$ are drawn next
to the corresponding nodes of standard Kostant modules.

\begin{figure}[ht]
\centering
\begin{pspicture}(-.3,-.35)(.2,15)
\psset{linewidth=.5pt,labelsep=8pt,nodesep=0pt}
\small
$
\psset{labelsep=7pt}
\cnode*(1.5,13){.07}{c1} 
\uput[d](1.5,12.96){a}
\cnode*(2.25,13){.07}{c2} 
\uput[d](2.25,13.05){b}
\cnode*(3,13){.07}{c3} 
\uput[d](3,12.96){c}
\cnode*(3.75,13){.07}{c4} 
\uput[d](3.75,13.05){d}
\ncline{c1}{c2}
\ncline{c3}{c4}
\psline(1.38,12.88)(1.62,13.12)
\psline(1.38,13.12)(1.62,12.88)
\psline(2.25,12.94)(3,12.94)
\psline(2.25,13.06)(3,13.06)
\psline(2.56,12.88)(2.68,13)
\psline(2.56,13.12)(2.68,13)
\cnode*(0,0){.07}{a1}
\pscircle(0,0){.2}
\cnode*(0,1){.07}{a2}
\pscircle(0,1){.2}
\cnode*(0,2){.07}{a3}
\pscircle(0,2){.2}
\cnode*(0,3){.07}{a4}
\pscircle(0,3){.2}
\cnode*(-.75,3.75){.07}{a5}
\pscircle(-.75,3.75){.2}
\cnode*(.75,3.75){.07}{a6}
\pscircle(.75,3.75){.2}
\cnode*(-1.5,4.5){.07}{a7}
\pscircle(-1.5,4.5){.2}
\cnode*(0,4.5){.07}{a8}
\cnode*(-.75,5.25){.07}{a9}
\cnode*(.75,5.25){.07}{a10}
\cnode*(0,6){.07}{a11}
\cnode*(1.5,6){.07}{a12}
\cnode*(0,7){.07}{a13}
\cnode*(1.5,7){.07}{a14}
\cnode*(-.75,7.75){.07}{a15}
\cnode*(.75,7.75){.07}{a16}
\cnode*(-1.5,8.5){.07}{a17}
\cnode*(0,8.5){.07}{a18}
\cnode*(-.75,9.25){.07}{a19}
\cnode*(.75,9.25){.07}{a20}
\cnode*(0,10){.07}{a21}
\cnode*(0,11){.07}{a22}
\cnode*(0,12){.07}{a23}
\cnode*(0,13){.07}{a24}
\pscircle(0,13){.2}
\ncline{a1}{a2}
\ncline{a2}{a3}
\ncline{a3}{a4}
\ncline{a4}{a5}
\ncline{a4}{a6}
\ncline{a5}{a7}
\ncline{a5}{a8}
\ncline{a6}{a8}
\ncline{a7}{a9}
\ncline{a8}{a9}
\ncline{a8}{a10}
\ncline{a9}{a11}
\ncline{a10}{a11}
\ncline{a10}{a12}
\ncline{a11}{a13}
\ncline{a11}{a14}
\ncline{a12}{a13}
\ncline{a12}{a14}
\ncline{a13}{a15}
\ncline{a14}{a16}
\ncline{a13}{a16}
\ncline{a15}{a17}
\ncline{a15}{a18}
\ncline{a16}{a18}
\ncline{a17}{a19}
\ncline{a18}{a19}
\ncline{a18}{a20}
\ncline{a19}{a21}
\ncline{a20}{a21}
\ncline{a21}{a22}
\ncline{a22}{a23}
\ncline{a23}{a24}
\uput[r](-0.15 ,.5){a}
\uput[r](-0.15 ,1.5){b}
\uput[r](-0.15 ,2.5){c}
\uput[r](-0.60,3.55){b}
\uput[l](0.60,3.55){d}
\uput[r](-1.35,4.25){a}
\uput[l](-0.15,4.30){d}
\uput[r](0.15,4.30){b}
\uput[l](-0.90,5.05){d}
\uput[r](-0.60,5.00){a}
\uput[l](0.60,5.00){c}
\uput[l](-0.15,5.75){c}
\uput[r](0.15,5.75){a}
\uput[l](1.35,5.80){b}
\uput[l](0.15,6.5){b}
\uput[r](1.35,6.5){a}
\uput[r](-0.60,7.50){c}
\uput[l](0.60,7.50){a}
\uput[r](0.90,7.55){b}
\uput[r](-1.35,8.30){d}
\uput[l](-0.15,8.25){a}
\uput[r](0.15,8.25){c}
\uput[l](-0.90,9.00){a}
\uput[r](-0.60,9.05){d}
\uput[l](0.60,9.05){b}
\uput[l](-0.15,9.80){b}
\uput[r](0.15,9.80){d}
\uput[r](-0.15 ,10.5){c}
\uput[r](-0.15 ,11.5){b}
\uput[r](-0.15 ,12.5){a}
\psset{labelsep=7pt}
\cnode*(-4.5,4.5){.07}{c1} 
\cnode*(-3.75,4.5){.07}{c2} 
\cnode*(-3,4.5){.07}{c3} 
\ncline{c1}{c2}
\psline(-4.62,4.38)(-4.38,4.62)
\psline(-4.62,4.62)(-4.38,4.38)
\psline(-3.75,4.44)(-3,4.44)
\psline(-3.75,4.56)(-3,4.56)
\psline(-3.44,4.38)(-3.32,4.5)
\psline(-3.44,4.62)(-3.32,4.5)
\psset{labelsep=7pt}
\cnode*(1.5,2){.07}{c1} 
\cnode*(2.25,2){.07}{c2} 
\ncline{c1}{c2}
\psline(1.38,1.88)(1.62,2.12)
\psline(1.38,2.12)(1.62,1.88)
\psset{labelsep=7pt}
\cnode*(1.5,1){.07}{c1} 
\psline(1.38,0.88)(1.62,1.12)
\psline(1.38,1.12)(1.62,0.88)
\rput(1.5,0){\varnothing}
$
\end{pspicture}
\caption{Kostant modules in a regular block for $(F_4,C_{3})$}
\label{F:F4C3}
\end{figure}
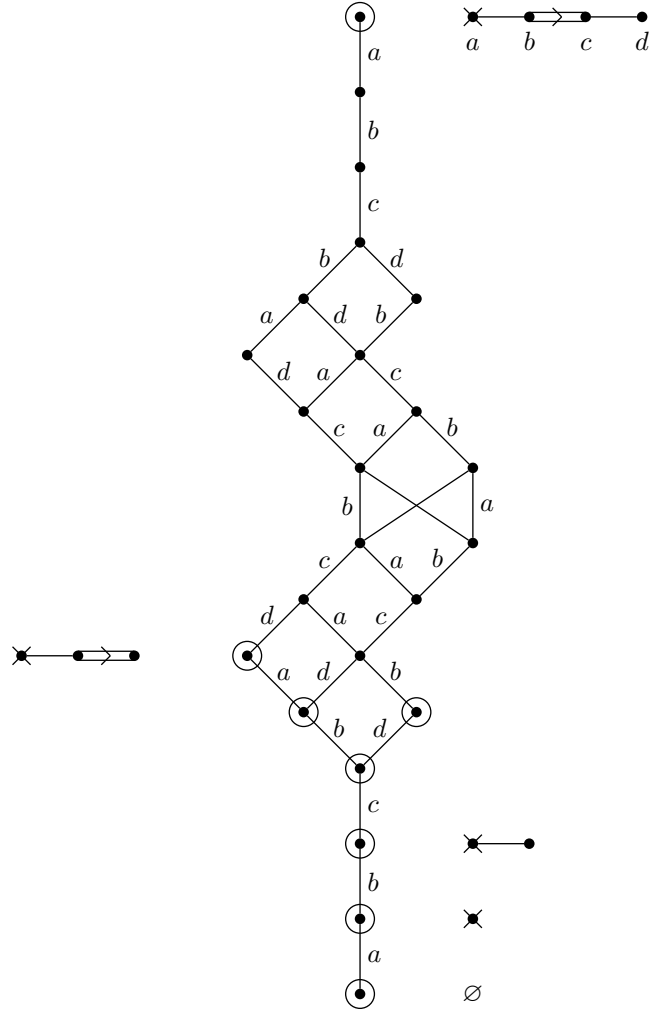
\end{ex}

\subsection{\it Enright-Shelton equivalences.} We conclude this section by giving an interpretation of standard Kostant modules via certain equivalences of categories that were introduced by Enright and Shelton in \cite{EnSh:87}. 
Fix $I\subset\Delta$, and let $\phi(I)$ be as in Section~\ref{SS:standardKostant}. Define the truncated category $\O_t(\fg,\fp,\phi(I))$ as the full subcategory of $\O(\fg,\fp,\reg)=\O_{S}^{\reg}$ consisting of all modules $V$ such that $[V: L_{w}]\not=0$ implies $w \leq \phi(I)$.
The simple modules in $\O_t(\fg,\fp,\phi(I))$ are then the modules $L_{w}$
with $w\leq  \phi(I)$. Let $\fl'=[\fl,\fl]$, where $\fl$ is the Levi 
subalgebra of $\fq=\fp_{I}$ as in Section \ref{SS:subdiagrams}.
In \cite{EnSh:87}, Enright and Shelton proved 
(see also \cite{Soe:89}) the following result. 

\begin{thm}\Label{P:ESequiv}
There is an equivalence of categories between $\O(\fl',\fl'\cap \fp,\reg)$ and  
$\O_{t}(\fg,\fp,\phi(I))$, sending the generalized Verma module $N_{w}'$ to $N_{w}$ and the simple module $L'_{w}$ to $L_{w}$ for all $w\in {}^{S'}W' \subset {}^{S}W$.
\end{thm}

Note that according to our definition of $\O(\fl',\fl'\cap \fp,\reg)$, the simple module in $\O(\fl',\fl'\cap \fp,\reg)$ corresponding to the maximal element of ${}^{S'}W'$ is finite dimensional.

\section{Regular Blocks for Maximal Parabolics in Simply Laced Type} 
\Label{S:ADEmaxreg}

\subsection{} In this section we assume that $\fg$ is simply-laced, and that $S^c=\{\a\}$ (i.e.\ that $\fp_S$ is a maximal parabolic subalgebra of $\fg$). We work in a fixed regular block $\O_S^{\text{reg}}$. Under these assumptions, we prove that there is a bijection between the set of connected subdiagrams of $\D$ containing $\a$ (recall our convention that we consider the empty subdiagram to belong to this set, since we really mean the set of subdiagrams with no $S$-trivial component) and the set of  Kostant modules in $\O_S^{\text{reg}}$.
We outline the idea of the proof. According to Proposition \ref{P:subdiagtoKostant}, we have an injection of the first set into the second. Recall that $L_w$ is a Kostant module (if and) only if the Poincar\'e polynomial $P_w(t)$ is palindromic. Set $\Pal=\{\, w\in\SW \mid P_w(t) \text{ is palindromic}\,\}$. 
We will show that 
\begin{equation} \Label{E:Palineq}
| \Pal | \le \#\{ \text{connected subdiagrams of } \D \text{ containing } \a \}.
\end{equation}
This will clearly prove that all three sets have the same cardinality.

\begin{thm} \Label{T:ADEmaxreg}
Let $\fg$ be a simply-laced simple Lie algebra, and let $S\subset\Delta$ with $S^c=\{\a\}$. Let $\D$ be the corresponding Dynkin diagram with one crossed node. Fix a regular block $\O_S^{\reg}$ of category $\O_S$. Then there is a
bijection
$$
\{\,\text{Kostant modules in }\O_S^{\text{reg}}\,\} \leftrightarrow 
\left\{
\begin{array}{c}
\text{connected subdiagrams}\\
\text{of } \D \text{ containing } \a\\  
\end{array}
\right\}
$$
(where we declare the empty subdiagram to belong to the set on the right).
\end{thm}

\begin{proof}
It suffices to prove \eqref{E:Palineq}. Certainly the identity element $e\in \Pal$. Fix $w\in\SW$ with $l(w)>0$. As $s_\a$ is the unique element of $\SW$ of length 1, and $s_\a \le w$, we have
\begin{equation} \Label{E:Pw}
P_w(t)=1+t+\dots+ct^{l(w)-1}+t^{l(w)} \quad \text{ for some } c\in\N.
\end{equation}
A necessary condition for $w\in\Pal$ is $c=1$. In particular, if $w\in\Pal$ then 
\begin{equation} \Label{E:onesimple}
\text{there is at most one simple reflection $s$ such that } w>ws\in\SW.
\end{equation}

In the classical types, we use Proctor's parameterization \cite{Proc:82} for $\SW$. (For reasons of personal preference we reverse his strings, left-to-right, and also reverse the Bruhat order, so that $e$ is the least element.) Assume $\a=\a_r$ in the Bourbaki ordering of simple roots.

First suppose we are in type $A_{n-1}$. Then $\SW$ consists of all permutations $w$ of $e=1^r 0^{n-r}$. The order-generating relations are given by
\begin{equation} \Label{E:Sij}
S_{ij}(w) < w \quad \text{ if } 1\le i<j\le n \text{ and } w_i<w_j,
\end{equation}
where $S_{ij}$ acts on the $n$-tuple $w$ by interchanging $w_i$ and $w_j$. Here $S_{ij}$ is the simple reflection $s_i$ when $j=i+1$. Evidently the elements satisfying \eqref{E:onesimple} are precisely those of the form
\begin{equation} \Label{E:AKostant}
w=1\dots \overset{k}{1}\, 0\dots0\, 1\dots \overset{l}{1}\, 0\dots 0,\quad 0\le k<r<l\le n \ \text{ or }\  k=r=l.
\end{equation}
On the other hand the connected subdiagrams of $\D$ containing $\a$ are obtained by ``pruning'' (deleting) nodes $k$ and $l$ from $\D$, where $0\le k< r<l\le n$ or $k=r=l$, and then selecting the connected component containing node $r$ in the remaining graph. (Here, the case $k=r=l$ produces the empty graph, which corresponds to the Kostant module $L_e$.) Since these two sets are parameterized by the identical values of $k$ and $l$, they have the same cardinality. (In fact, it is easy to see that the element in \eqref{E:AKostant} is the maximal element of the subposet of $\SW$ corresponding to the connected sub-Dynkin diagram obtained by pruning nodes $k$ and $l$ as described.) This proves the desired bijection in type $A$.

Next suppose we are in type $D_n$. Assume first that $r< n-1$. Then $\SW$ consists of all signed permutations of $e=1^r 0^{n-r}$. We denote $-1$ by $\bar 1$ when depicting signed permutations. The order-generating relations are given by $S_{ij}$ as in \eqref{E:Sij}, together with
\begin{equation} \Label{E:SNij}
SN_{ij}(w) < w \quad \text{ if } 1\le i<j\le n \text{ and } w_i+w_j<0,
\end{equation}
where $SN_{ij}$ acts by interchanging and negating $w_i$ and $w_j$. Here $SN_{n-1,n}$ is the simple reflection $s_n$. The elements satisfying \eqref{E:onesimple} are of the following types:
\begin{equation} \Label{E:DKostant}
\begin{aligned}
w^1&=1\dots \overset{k}{1}\, 0\dots0\, 1\dots \overset{l}{1}\, 0\dots 0,\quad 0\le k<r<l\le n \ \text{ or }\  k=r=l,\\
w^2&=1\dots \overset{k}{1}\, 0\dots 0\, 1\dots 1 \bar 1,\quad 0\le k<r,\\
w^3&=1\dots \overset{k}{1}\, \bar 1 \dots \bar 1 \, 0\dots 0,\quad 0\le k<r,\\
w^4&=1\dots \overset{k}{1}\, 0\dots \overset{l}{0}\, \bar 1 \dots \overset{m}{\bar 1}\, 1\dots \overset{p}{1} \, 0\dots 0,\quad 0\le k\le l<m<p\le n\\
w^5&=1\dots \overset{k}{1}\, 0\dots 0\, \overset{l}{\bar 1} \dots \bar1 \bar1,\quad 0\le k \le r-2,\\
w^6&=1\dots \overset{k}{1}\, 0\dots 0\, \bar1 \dots \overset{l}{\bar 1}\, 0\dots 0,\quad 0\le k<r<l\le n-2.
\end{aligned}
\end{equation}
The element $w^1$ is the maximal element of the sub-poset of $\SW$ corresponding to the subdiagram of $\D$ consisting of nodes $k+1$ through $l-1$ (or the empty subdiagram, if $k=r=l$). The element $w^2$ corresponds to the subdiagram consisting of nodes $k+1$ through $n-2$ and $n$. The element $w^3$ corresponds to the subdiagram consisting of nodes $k+1$ through $n$. These exhaust all the possible connected subdiagrams of $\D$ containing node $r$. Thus it remains to show that elements of the form $w^4, w^5, w^6$ do not parameterize Kostant modules. That is, we must show that their Poincar\'e polynomials are not palindromic.

Let $w=w^4$. Then $ws_m\rightarrow w$. But also, by \cite[Cor.\ 5D]{Proc:82}, $S_{l,m+1}(w)\rightarrow w$ if $k<l$, and $S_{m,p+1}(w)\rightarrow w$ if $p<n$. At least one of these conditions holds since $r<n-1$. Thus, in the notation of \eqref{E:Pw}, $c\ge 2$ and $P_w(t)$ is not palindromic.

Let $w=w^5$. Note that the condition on $k$ implies $l\le n-1$; the case $l=n$ is of type $w^2$. Then $ws_n\rightarrow w$. But also, by \cite[Cor.\ 5D]{Proc:82}, $SN_{l-1,n}(w)\rightarrow w$. Again, $P_w(t)$ is not palindromic.

Let $w=w^6$. This is the most challenging case, because we must go several levels down from $w$ (and up from $e$) to detect the failure of palindromicity. As a warm-up, consider the case $r=1$, where $w=0\dots 0\,\overset{l}{\bar 1}\, 0\dots 0$ with $1<l\le n-2$. Write $[e,w]_i = \{\,x\le w\mid l(x)=i\,\}$. We claim that $\#[e,w]_{l(w)-i}=\#[e,w]_i=1$ for $0\le i<n-l$ but $\#[e,w]_{l(w)-(n-l)}=2$ while $\#[e,w]_{n-l}=1$, and thus $P_w(t)$ is not palindromic. Indeed, for $i\le n-l$, $[e,w]_i = \{\,0\dots 0\,\overset{i+1}{1}\, 0\dots 0\,\}$ and for $i<n-l$, $[e,w]_{l(w)-i}=\{\,0\dots 0\,\overset{l+i}{\bar 1}\, 0\dots 0\,\}$, but $[e,w]_{l(w)-(n-l)}=\{\,0\dots 0\,\bar 1\,\}\cup\{\, 0\dots 0\,1\,\}$. Returning to general $r<n-1$, we claim that $\#[e,w]_{l(w)-i}\ge\#[e,w]_i$ for $0\le i<n-l$ but $\#[e,w]_{l(w)-(n-l)}>\#[e,w]_{n-l}$. In fact, any element $x$ of $[e,w]_i$ for $0\le i\le n-l$ is obtained from $e$ by moving some of the last $r-k$ 1's to the right a total of $i$ positions. (Moving any of the first $k$ 1's would produce an element not dominated by $w$. And there are no non-simple coverings among the elements of length at most $n-l$.) To any such $x$ there is an element $x'$ of $[e,w]_{l(w)-i}$ obtained by moving the corresponding $\bar 1$'s in $w$ the same number of positions to the right. (The last 1 in $e$ corresponds to the last $\bar 1$ in $w$, and so on.) The map $x\mapsto x'$ is clearly an injection. On the other hand, the element $1\dots \overset{k}{1}\, 0\dots 0\, \bar1 \dots \overset{l-1}{\bar 1}\, 0\dots 0\, 1 \in [e,w]_{l(w)-(n-l)}$ is not in the image of this map. This shows that $P_w(t)$ is not palindromic, and completes the proof for type $D_n$ when $r<n-1$.

The labeled posets $\SW$ for $r=n-1$ and $r=n$ are isomorphic, so it suffices to consider the case $r=n$. Then $\SW$ consists of all signed permutations of $e=1^n$ containing an even number of $\bar 1$'s. The elements satisfying \eqref{E:onesimple} are of the following types:
\begin{equation} 
\begin{aligned}
w^1&=1\dots 1\\
w^2&=1\dots \overset{k}{1}\, \bar 1\, 1\dots 1\, \bar 1,\quad 0\le k\le n-2,\\
w^3&=1\dots \overset{k}{1}\, \bar 1 \dots \bar 1 \pm\! 1,\quad 0\le k< n-2,\\
w^4&=1\dots \overset{k}{1}\,  \bar 1 \dots \overset{m}{\bar 1}\, 1\dots 1 \pm\!1 ,\quad 0\le k<m-1,\ m< n-1,
\end{aligned}
\end{equation}
with $\pm$ chosen in $w^3$ and $w^4$ to make the total number of $-1$'s even. Then $w^1$ corresponds to the empty subdiagram of $\D$, $w^2$ to the subdiagram containing nodes $\{k+1,\dots,n\}\smallsetminus\{n-1\}$, and $w^3$ to the subdiagram containing nodes $\{k+1,\dots,n\}$. This accounts for all connected subdiagrams of $\D$ containing node $n$, so we must show that $w=w^4$ has non-palindromic Poincar\'e polynomial. Notice that $\SW$ contains a unique element of length 2, namely $1\dots 1 \bar 1 1 \bar 1$. However, there are two elements in $[e,w]_{l(w)-2}$, namely $ws_ms_{m-1}$ and $ws_ms_{m+1}$ (or $ws_ms_{m+2}$ if $w=1\dots 1 \bar 1 \dots \bar 1  1 \bar1$). Thus the coefficients of $t^2$ and $t^{l(w)-2}$ differ in $P_w(t)$.

Finally assume we are in type $E_n$. We wrote a computer program to calculate $\SW$ (using the characterization $\SW=\{\,w\in W\mid w^{-1}(\Phi_S^+)\subset \Phi^+\,\}$) and then find the elements in $\SW$ having palindromic Poincar\'e polynomial. (More complicated computer analyses of rational smoothness for parabolic Schubert varieties in type $E$ have been done by Billey-Postnikov \cite{BiPo:05}.) The Bruhat order was determined via the following recursive algorithm.  Given $x, w\in\SW$ with $w\ne e$, find a simple reflection $s$ such that $w>ws=:w'$. Then necessarily $w'\in\SW$, by \cite{Deo:77}. Set $x'$ equal to $x$ or $xs$, whichever is smaller. Then $x\le w$ if and only if $x'\le w'$, by op.\ cit. It is easy to calculate by hand the number of connected subdiagrams of $D$ containing node $r$, for each choice of $r$. In every case, we found the number of elements with palindromic Poincar\'e polynomial to be the same as the number of such connected subdiagrams.  These numbers are listed in Table \ref{T:Emax}.
\begin{table}[h]
\begin{center}
\begin{tabular}{cccc}
$r$ &  $E_6$ & $E_7$ & $E_8$ \\[2pt] \hline 
1 & 9   & 11 & 13  \\
2 & 11 & 14 & 17  \\
3 & 15 & 19 & 23  \\
4 & 19 & 25 & 31  \\
5 & 15 & 22 & 29  \\
6 & 9   & 17 & 25  \\
7 &      & 10 & 19  \\
8 &      &       & 11  \\ \hline
\end{tabular}
\end{center}
\caption{Number of Kostant modules for maximal parabolics} \Label{T:Emax}
\end{table}
\end{proof}

\section{Arbitrary Blocks for Hermitian Symmetric Pairs} 
\Label{S:HSsing}

\subsection{} In this section we assume that $(\Phi,\Phi_S)$ corresponds to a Hermitian symmetric pair.  We prove that in any (nonempty) integral block $\O_S^\mu$, the Kostant modules are parameterized by a set of the form $\{$connected subdiagrams of $\D'$ containing $\a' \}$, for a certain Dynkin diagram $\D'$ with a single crossed node $\a'$. The key tool is an equivalence of categories due to Enright (based on work of Enright and Shelton). In section \ref{SS:unified} we present a unified construction of the diagram $\D'$.

As remarked in the introduction, the Kostant modules in regular (integral) Hermitian symmetric categories were first worked out by Collingwood \cite{Col:85}. We give new proofs of his results.

\subsection{\it The regular case} Let $(\Phi,\Phi_S)$ correspond to a Hermitian symmetric pair.  Then $S$ is a maximal proper subset of $\Delta$; put $S^c=\{\a\}$ as before. Assume in this subsection that $\mu$ is regular.

\begin{figure}[ht]
\centering
\begin{pspicture}(-.3,-.35)(2.3,4)
$
\cnode*(-5,2){.07}{a}
\cnode*(-4,2){.07}{b}
\pnode(-3.5,2){c}
\pnode(-3,2){d}
\pnode(-2.5,2){w}
\pnode(-2,2){x}
\cnode*(-1.5,2){.07}{y}
\cnode*(-.5,2.5){.07}{z}
\cnode*(-5,3){.07}{a1}
\cnode*(-4,3){.07}{b1}
\pnode(-3.5,3){c1}
\pnode(-3,3){d1}
\pnode(-2.5,3){w1}
\pnode(-2,3){x1}
\cnode*(-1.5,3){.07}{y1}
{\psset{linewidth=.5pt,labelsep=8pt}
\ncline{a}{b}
\ncline{b}{c}
\ncline{x}{y}
\ncline{y}{z}
\ncline{a1}{b1}
\ncline{b1}{c1}
\ncline{x1}{y1}
\ncline{y1}{z}
\psline(-5.15,2.15)(-4.85,1.85)
\psline(-5.15,1.85)(-4.85,2.15)
{\psset{linewidth=2pt,linestyle=dotted,dotsep=6pt}
\ncline{d}{w}
\ncline{d1}{w1}
}
{\small
\uput[u](-5,3){2n-1}
\uput[d](-.5,2.5){n}
\uput[d](-5,2){1}
}
}
\psline{->}(-2.75,1.5)(-2.75,.5)
\cnode*(-5,0){.07}{a2}
\cnode*(-4,0){.07}{b2}
\pnode(-3.5,0){c2}
\pnode(-3,0){d2}
\pnode(-2.5,0){w2}
\pnode(-2,0){x2}
\cnode*(-1.5,0){.07}{y2}
\cnode*(-.5,0){.07}{z2}
{\psset{linewidth=.5pt,labelsep=8pt}
\ncline{a2}{b2}
\ncline{b2}{c2}
\ncline{x2}{y2}
{\psset{doubleline=true,doublesep=3pt}
\ncline{-}{y2}{z2}
}
\psline(-.9,.15)(-1.1,0)(-.9,-.15)
\psline(-5.15,.15)(-4.85,-.15)
\psline(-5.15,-.15)(-4.85,.15)
{\psset{linewidth=2pt,linestyle=dotted,dotsep=6pt}
\ncline{d2}{w2}
}
{\small
\uput[d](-5,0){1}
\uput[d](-.5,0){n}
}
}
\cnode*(2,2.5){.07}{a3}
\cnode*(3,2.5){.07}{b3}
\pnode(3.5,2.5){c3}
\pnode(4,2.5){d3}
\pnode(4.5,2.5){w3}
\pnode(5,2.5){x3}
\cnode*(5.5,2.5){.07}{y3}
\cnode*(6.5,2){.07}{z3}
\cnode*(6.5,3){.07}{zz3}
{\psset{linewidth=.5pt,labelsep=8pt}
\ncline{a3}{b3}
\ncline{b3}{c3}
\ncline{x3}{y3}
\ncline{y3}{z3}
\ncline{y3}{zz3}
\psline(6.65,2.15)(6.35,1.85)
\psline(6.65,1.85)(6.35,2.15)
{\psset{linewidth=2pt,linestyle=dotted,dotsep=6pt}
\ncline{d3}{w3}
}
{\small
\uput[d](2,2.5){1}
\uput[u](6.5,3){n}
\uput[d](6.5,2){n+1}
}
}
\psline{->}(4.25,1.5)(4.25,.5)
\cnode*(2,0){.07}{a4}
\cnode*(3,0){.07}{b4}
\pnode(3.5,0){c4}
\pnode(4,0){d4}
\pnode(4.5,0){w4}
\pnode(5,0){x4}
\cnode*(5.5,0){.07}{y4}
\cnode*(6.5,0){.07}{z4}
{\psset{linewidth=.5pt,labelsep=8pt}
\ncline{a4}{b4}
\ncline{b4}{c4}
\ncline{x4}{y4}
{\psset{doubleline=true,doublesep=3pt}
\ncline{-}{y4}{z4}
}
\psline(5.9,.15)(6.1,0)(5.9,-.15)
\psline(6.35,.15)(6.65,-.15)
\psline(6.35,-.15)(6.65,.15)
{\psset{linewidth=2pt,linestyle=dotted,dotsep=6pt}
\ncline{d4}{w4}
}
{\small
\uput[d](2,0){1}
\uput[d](6.5,0){n}
}
}
$
\end{pspicture}
\caption{Simply-laced covers}
\label{F:fold}
\end{figure}
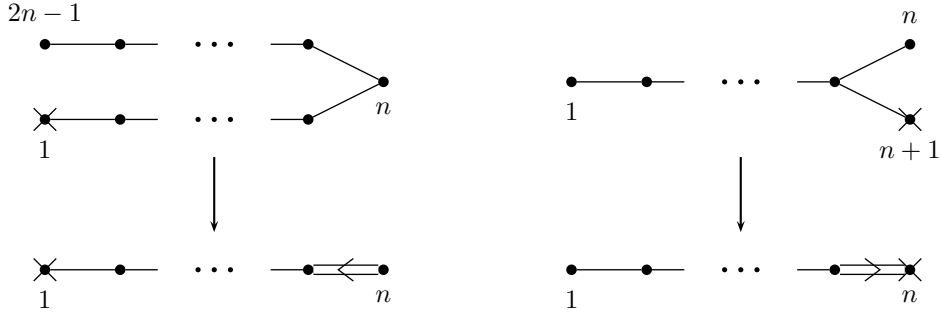

We need some additional notation in case $\D$ is not simply laced; i.e., when $\D=(B_n,B_{n-1})$ (resp.\ $(C_n,A_{n-1})$). Let $\D^\vee$ be the Dynkin diagram dual to $\D$, and let $\a^\vee$ be the simple root in $\D^\vee$ dual to $\a$. (The fundamental dominant weight associated to $\a$ is cominuscule, while that associated to $\a^\vee$ is minuscule.) Let $\Phi^\vee$ be the root system of $\D^\vee$ with simple roots $\Delta^\vee$. Let $ S^\vee=\Delta^\vee \smallsetminus \{\a^\vee \}$. Define the \emph{simply-laced cover} of $\D^\vee$ to be the Dynkin diagram $\widetilde {\D^\vee}$ of type $A_{2n-1}$ (resp.\ $D_{n+1}$). It can be useful to imagine $\widetilde {\D^\vee}$ ``folding'' at node $n$ (resp.\ $n-1$) to ``cover'' $\D^\vee$, as depicted in Figure \ref{F:fold}. Let $\widetilde{\Phi^\vee}$ be the root system of $\widetilde {\D^\vee}$ with simple roots $\widetilde{\Delta^\vee}$. Let $\widetilde {S^\vee}=\widetilde{\Delta^\vee} \smallsetminus \{\widetilde{\a^\vee} \}$, where $\widetilde{\a^\vee}=\a_1$ (resp.\ $\a_{n+1}$). Thus we have associated to $\D=(B_n,B_{n-1})$ (resp.\ $(C_n,A_{n-1})$) a  simply-laced Hermitian symmetric pair  $\widetilde{\D^\vee} = (A_{2n-1},A_{2n-2})$ (resp.\ $(D_{n+1},A_n)$).

\begin{thm} \Label{T:HSreg}
Let $\O_S^{\text{reg}}$ be a regular block for a Hermitian symmetric pair. Then there is a Dynkin diagram $\D'$ with one crossed node $\a'$, such that 
$$
\{\,\text{Kostant modules in }\O_S^{\text{reg}}\,\} \leftrightarrow 
\left\{
\begin{array}{c}
\text{connected subdiagrams}\\
\text{of } \D' \text{ containing } \a'\\  
\end{array}
\right\}
$$
(where we declare the empty subdiagram to belong to the set on the right).
If $\D$ is simply-laced, then $\D'=\D$ and $\a'=\a$; otherwise, $\D'=\widetilde {\D^\vee}$ and $\a'=\widetilde{\a^\vee}$.
\end{thm}

\begin{proof}
By Theorem \ref{T:ADEmaxreg}, if $\D$ is simply-laced, there is nothing to prove. If not, let $\D^\vee,\ \widetilde {\D^\vee}$, etc.\ be as above. Since the Weyl groups for $B_n$ and $C_n$ are identical, so are the labeled posets $\SW$ associated to $\D$ and $\D^\vee$ (with crossed nodes $\a$ and $\a^\vee$, respectively), and they have the same KL polynomials and hence the same parameters for the Kostant modules. Thus we may replace $\D$ with $\D^\vee$. Now let $G^\vee$ be a connected simple algebraic group over $\C$ having root system $\Phi^\vee$, $P^\vee\subset G^\vee$ a parabolic subgroup corresponding to $S^\vee$, and similarly for $\widetilde {G^\vee}$ and $\widetilde {P^\vee}$. The flag varieties $G^\vee/P^\vee$ and $\widetilde {G^\vee}/\widetilde {P^\vee}$ are isomorphic, via an isomorphism (coming from the embeddings $Sp(2n,\C)\hookrightarrow SL(2n,\C)$ and $SO(2n+1,\C)\hookrightarrow SO(2n+2,\C)$) which identifies Schubert varieties in the two spaces \cite[Sec.\ 3.1]{BrPo:99}.  Thus, by our identification of Kostant modules with rationally smooth Schubert varieties, and Theorem \ref{T:ADEmaxreg}, we have that the Kostant modules are parameterized by connected subdiagrams of $\widetilde {\D^\vee}$ containing $\widetilde{\a^\vee}$ (together with the empty diagram). 
\end{proof}

\subsection{\it The singular case}  Maintain the assumptions on $S$ from the previous subsection, but now assume that $\mu$ is singular. Define $J$ as in Section \ref{SS:posets}. If $J$ contains two simple roots adjacent in $\D$, then $\O_S^\mu$ is empty, by \cite[Cor.~4.2]{BoNa:05} or \cite[Lem.~3.1a]{Enr:88}. So we may assume that $\Phi_J$ is of type $A_1\times\dots\times A_1$ ($t$ factors, for some $t\in\N$). Suppose we are in one of the three ``complicated'' Hermitian symmetric cases, $(A_n,A_{r-1}\times A_{n-r}), (C_n,A_{n-1})$, or $(D_n,A_{n-1})$. In 1987 Enright and Shelton \cite{EnSh:87} proved in the semi-regular case ($t=1$) that the category $\O(\Phi,\Phi_S,\Phi_J)$ is  equivalent to a regular block of a Hermitian symmetric category of rank $n-2$ (or, in one case, a direct sum of two copies of such a block). The following year Enright \cite{Enr:88} realized that essentially the same proof handles the case of general $J$, with the rank of the regular block being $n-2t$. As mentioned above, these equivalences are key to understanding the Kostant modules in these three cases.

To get a feel for how the equivalences work, it may be helpful to look at an example. Consider type $A_{n-1}$ with $S^c=\{r\}$ and $J=\{j_1,\dots,j_t\}$ (viewing $S, J \subset \{1,\dots,n-1\}$ in the obvious way). After applying a translation functor, and shifting by a constant vector to make all entries non-negative, we may assume that $\mu+\rho=(a_1,\dots,a_n)$ (in the $\e_i$-coordinates) with $a_1=0$, $a_{i+1}=a_i+1$ if $i\notin J$, and $a_{i+1}=a_i$ if $i\in J$. The permutations of $\mu+\rho$ which are highest weights plus $\rho$ of simple modules in $\O(\Phi,\Phi_S,\Phi_J)$ have their first $r$ entries strictly decreasing, and likewise their last $n-r$ entries strictly decreasing. This implies that, of each matched pair of equal entries $a_j=a_{j+1}\ (j\in J)$, one entry must occur in the first $r$ coordinates and the other in the last $n-r$. Now these matched pairs can all be deleted; this operation is clearly invertible. The resulting $(n-2t)$-tuples correspond to the weights for a regular block of the Hermitian symmetric category $\O(A_{n-1-2t},A_{r-1-t}\times A_{n-1-r-t})$. Enright's result is that this bijection on weights arises from an actual equivalence of categories. Similar weight analyses can be easily made for the other two ``complicated'' cases.

\begin{thm} \Label{T:HSsing} Let $(\Phi,\Phi_S)$ correspond to a Hermitian symmetric pair. In case $(\Phi,\Phi_S)=(C_n,A_{n-1})$, assume $J$ does not contain the long simple root. Then whenever $\O(\Phi,\Phi_S,\linebreak[0]\Phi_J)$ is a nonempty singular block, there is a Dynkin diagram $\D'$ with a single crossed node $\a'$ such that the set of Kostant modules in $\O(\Phi,\Phi_S,\Phi_J)$ is in bijection with one or two copies of $\{$connected subdiagrams of $\D'$ containing $\a'\}$. The pairs $(\D', \a')$, and the number of copies are given explicitly in Table~\ref{T:HSsingdata}. (The case $\D'=\varnothing$ corresponds to the block $\O(\varnothing,\varnothing)$ having only one simple module.)

\end{thm}

\begin{table}[h] 
\begin{tabular}{cccccc}
$\D$ &  $\a$ & $|J|$ & $\D'$ & $\a'$ & \# Copies \\[2pt] \hline 
$(A_n,A_{r-1}\times A_{n-r})$ & $\a_r$   & $t$ & $(A_{n-2t},A_{r-t-1}\times A_{n-r-t})$ & $\a_{r-t}$ & 1\\
$(B_n,B_{n-1})$ & $\a_1$   & 1 (short) & $\varnothing$  & -- & 1\\
$(B_n,B_{n-1})$ & $\a_1$   & 1 (long) & $\varnothing$ & -- & 2\\
$(C_n,A_{n-1})$ & $\a_n$   & $t$ (all short) & $(D_{n+1-2t},A_{n-2t})$  & $\a_{n+1-2t}$ & 1 \\
$(C_n,A_{n-1})$ & $\a_n$   & $t$ (1 long) & $(D_{n+1-2t},A_{n-2t})$  & $\a_{n+1-2t}$ & 2 \\
$(D_n,D_{n-1})$ & $\a_1$   & 1 & $(A_1,\varnothing)$  & $\a_1$ & 1 \\
$(D_n,D_{n-1})$ & $\a_1$   & 2 ($\{n-1,n\}$) & $\varnothing$  & -- & 1 \\
$(D_n,A_{n-1})$ & $\a_n$   & $t$ & $(D_{n-2t},A_{n-2t-1})$  & $\a_{n-2t}$ & 1 \\
$(E_6,D_5)$ & $\a_6$   & 1 & $(A_5,A_4)$  & $\a_5$ & 1 \\
$(E_6,D_5)$ & $\a_6$   & 2 & $\varnothing$  & -- & 1 \\
$(E_7,E_6)$ & $\a_7$   & 1 & $(D_6,D_5)$  & $\a_1$ & 1 \\
$(E_7,E_6)$ & $\a_7$   & 2 & $(A_1,\varnothing)$  & $\a_1$ & 1 \\
$(E_7,E_6)$ & $\a_7$   & 3 ($\{2,5,7\}$) & $\varnothing$  & -- & 1 \\[2pt] \hline
\end{tabular}
\medskip
\caption{Data for singular Hermitian symmetric categories} \Label{T:HSsingdata}
\end{table}

\begin{rem}
The theorem should ``morally'' be true for the excluded case  $\O(C_n,A_{n-1},\Phi_J)$ where $J$ contains the long simple root. When such a category is nonempty, it splits into a direct sum of two blocks, each equivalent to $\O(C_{n-2|J|},A_{n-2|J|-1|},\varnothing)$, by \cite[Prop.~3.2(b)]{Enr:88}. So a simple module $L_w$ in the original category ought to be a Kostant module if it corresponds to a Kostant module in either of the direct summands. Unfortunately, the relevant interval $[v,w]$ in $\SWJ$ will usually include parameters $x$ from the other direct summand, which do not appear in the cohomology of $L_w$. Thus $L_w$ cannot satisfy Definition \ref{D:Kostant} of  Kostant modules. In Section \ref{SS:Extorder} we propose a new ordering on $\SWJ$ which will rectify this problem and make $L_w$ a Kostant module; see Section \ref{SS:HSsplit} for details. In anticipation of that discussion, we have included the data for this case in Table~\ref{T:HSsingdata}.
\end{rem}

\begin{proof}
First assume we are in type $(A_n,A_{r-1}\times A_{n-r}),\ (C_n,A_{n-1})$, or $(D_n,A_{n-1})$.
Then by \cite[Prop.\ 3.2(a)]{Enr:88} there is an equivalence of categories $\O(\Phi,\Phi_S,\Phi_J) \simeq \O(\Phi',\Phi'_{S'},\varnothing)$, where $(\Phi',\Phi'_{S'})$ corresponds to a Hermitian symmetric pair of the same type but of rank $n-2|J|$. The equivalence induces a bijection of partially ordered sets $\SWJ \leftrightarrow {}^{S'} W'$, and identifies correspondingly-parameterized simple modules, GVMs, etc. In particular (using the characterization in terms of palindromic Poincar\'e polynomials), the equivalence matches up Kostant modules in the two categories.  

Next assume we are in one of the remaining classical types, $(B_n,B_{n-1})$ or $(D_n,D_{n-1})$. The structure of the semi-regular categories ($J=\{j\}$) was described in \cite[Sec.~4.1]{BoNa:05}. In type $B$ if $j<n$ there are two simple modules in two separate blocks, while if $j=n$ there is just one simple module. In type $D$ for any $j$ there are two simple modules having an extension between them. Since the poset $\SW$ in type $B$ is a chain, the set $\SWJ$ is empty whenever $|J|\ge 2$, so the corresponding categories are all empty. In type $D$, $\SWJ$ is a ``chain with a diamond in the middle,'' the edges up from the bottom of the diamond being labeled by $n-1$ and $n$. So for $|J|\ge 2$ the category is empty unless $J=\{n-1,n\}$, in which case it contains a single simple module. The theorem follows from these analyses.

Finally assume we are in one of the exceptional cases, $(E_6,D_5)$ or $(E_7,E_6)$. The structure of the singular blocks of these categories was analyzed in \cite{EnSh:89}. Begin with $(E_6,D_5)$. The semi-regular categories are all equivalent, and have the same structure as the regular category for $(A_5,A_4)$: the poset of highest weights is a chain of length six, and if these are numbered 1 through 6 with $V_1$ the irreducible GVM, then $\dim\Ext^k(V_i,L_j)=\d_{k,j-i}$ \cite[Prop.~2.3]{EnSh:89}. In particular, every simple module is a Kostant module, as for $\O(A_5,A_4,\varnothing)$. If $J$ consists of any two orthogonal simple roots, it is easy to see by inspection of \cite[Fig.~3.1]{EnSh:89} that $\SWJ$ contains only a single node, so the block consists of a single simple module. Now consider $(E_7,E_6)$. The semi-regular categories are again all equivalent, and have the same Hasse diagram and $\Ext(V_x,L_w)$-structure as $\O(D_6,D_5,\varnothing)$ \cite[Prop.~3.9]{EnSh:89}. In particular, the Kostant modules are the same as for the regular $(D_6,D_5)$ category. Similarly the categories in which $J$ consists of two orthogonal simple roots are all equivalent, and have the same structure as $\O(A_1,\varnothing,\varnothing)$ \cite[Lemma~3.5]{EnSh:89}. Lastly, if $J$ consists of three orthogonal simples, then by inspection of \cite[Fig.~3.1]{EnSh:89}, the category is nonempty only if $J=\{2,5,7\}$ (in the Bourbaki labeling), in which case it contains only a single simple module.
\end{proof}

\subsection{\it A unified approach} \Label{SS:unified} In this subsection we present a unified construction of the diagram $\D'$ associated to a Hermitian symmetric diagram $\D$ with a single crossed node $\a$ in Theorem \ref{T:HSsing}. We follow the Enright-Shelton convention that if all roots are the same length, then they are all short.

Let $\{\g_1,\dots,\g_u\}$ be a maximal set of strongly orthogonal short roots in $\Phi(\fu)$, defined as follows. Let $\g_1$ be the highest short root in $\Phi(\fu)$, and for $i>1$ let $\g_i$ be the highest short root in $\Phi(\fu)$ that is orthogonal to $\g_1,\dots,\g_{i-1}$. (These are strongly orthogonal because $\fu$ is abelian, so the sum of two roots in $\Phi(\fu)$ is never a root.) For $0\le t\le u$ define $\Phi^{(t)}=\{\,\a\in\Phi\mid (\a,\g_i)=0 \text{ for } 1\le i\le t\,\}$. We construct a sequence of Dynkin diagrams $\D^{(t)}$ associated to the root systems $\Phi^{(t)}$ as follows. Starting with $\D^{(0)}=\D$, create an extended Dynkin diagram by attaching an extra node corresponding to $-\g_1$. Then delete $-\g_1$ and all nodes adjacent to it, to obtain $\D^{(1)}$. In general, $\D^{(t)}$ is obtained from $\D^{(t-1)}$ by the same process (attaching $-\g_t$). One can check that $\D^{(t)}$ again corresponds to a Hermitian symmetric pair. (The Dynkin diagram $\D^{(t)}$ may be disconnected, so that $\O(\Phi^{(t)}, \Phi^{(t)}\cap \Phi_S)$ is a direct product of the categories corresponding to its connected components, all but one of which will be trivial. We may therefore delete these $S$-trivial components, replacing $\D^{(t)}$ with its connected component containing $\a$. In two cases, $(A_n,A_{n-1})$ and $(B_n,B_{n-1})$, the simple root $\a$ is adjacent to $-\g_1$, so there is no component of $\D^{(1)}$ containing $\a$. In this situation we replace $\D^{(1)}$ with $\varnothing$, with the same meaning as explained in the statement of Theorem \ref{T:HSsing}. We do the same thing in general if $\a\notin \D^{(t)}$.) In fact, comparison of the (connected) diagrams so obtained with those in Table \ref{T:HSsingdata} yields the following unified statement.

\begin{thm} \Label{T:HSsingunif} Let $(\Phi,\Phi_S)$ correspond to a Hermitian symmetric pair. If $\O(\Phi,\Phi_S,\Phi_J)$ is a nonempty singular block with 
$|J|=t$, then the diagram $\D'$ of Theorem \ref{T:HSsing} is $\D^{(t)}$ if $\D^{(t)}$ is simply laced, otherwise the simply-laced cover of its dual.
\end{thm}

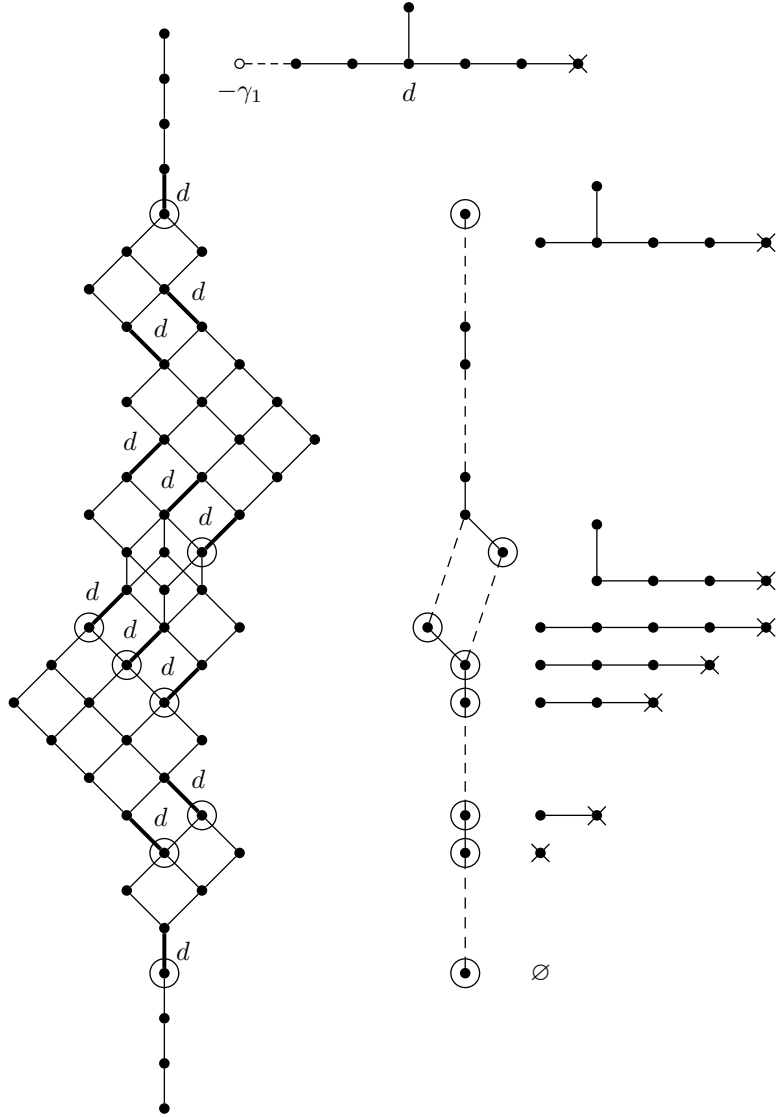
\begin{figure}[ht]
\centering
\begin{pspicture}(-.3,0)(1,15)
\psset{linewidth=.5pt,labelsep=3pt,nodesep=0pt}
\small
$
\cnode*(-2,.1){.07}{a1}
\cnode*(-2,.7){.07}{a2}
\cnode*(-2,1.3){.07}{a3} 
\cnode*(-2,1.9){.07}{a4} 
\pscircle(-2,1.9){.2}
\cnode*(-2,2.5){.07}{a5} 
\cnode*(-2.5,3){.07}{a6} 
\cnode*(-1.5,3){.07}{a7} 
\cnode*(-2,3.5){.07}{a8} 
\pscircle(-2,3.5){.2}
\cnode*(-1,3.5){.07}{a9} 
\cnode*(-2.5,4){.07}{a10} 
\cnode*(-1.5,4){.07}{a11} 
\pscircle(-1.5,4){.2}
\cnode*(-3,4.5){.07}{a12} 
\cnode*(-2,4.5){.07}{a13} 
\cnode*(-3.5,5){.07}{a14} 
\cnode*(-2.5,5){.07}{a15} 
\cnode*(-1.5,5){.07}{a16} 
\cnode*(-4,5.5){.07}{a17} 
\cnode*(-3,5.5){.07}{a18} 
\cnode*(-2,5.5){.07}{a19} 
\pscircle(-2,5.5){.2}
\cnode*(-3.5,6){.07}{a20} 
\cnode*(-2.5,6){.07}{a21} 
\pscircle(-2.5,6){.2}
\cnode*(-1.5,6){.07}{a22} 
\cnode*(-3,6.5){.07}{a23} 
\pscircle(-3,6.5){.2}
\cnode*(-2,6.5){.07}{a24} 
\cnode*(-1,6.5){.07}{a25} 
\cnode*(-2.5,7){.07}{a26} 
\cnode*(-2,7){.07}{a27} 
\cnode*(-1.5,7){.07}{a28} 
\cnode*(-2.5,7.5){.07}{a29} 
\cnode*(-2,7.5){.07}{a30} 
\cnode*(-1.5,7.5){.07}{a31} 
\pscircle(-1.5,7.5){.2}
\cnode*(-3,8){.07}{a32} 
\cnode*(-2,8){.07}{a33} 
\cnode*(-1,8){.07}{a34} 
\cnode*(-2.5,8.5){.07}{a35} 
\cnode*(-1.5,8.5){.07}{a36} 
\cnode*(-.5,8.5){.07}{a37} 
\cnode*(-2,9){.07}{a38} 
\cnode*(-1,9){.07}{a39} 
\cnode*(0,9){.07}{a40} 
\cnode*(-2.5,9.5){.07}{a41} 
\cnode*(-1.5,9.5){.07}{a42} 
\cnode*(-.5,9.5){.07}{a43} 
\cnode*(-2,10){.07}{a44} 
\cnode*(-1,10){.07}{a45} 
\cnode*(-2.5,10.5){.07}{a46} 
\cnode*(-1.5,10.5){.07}{a47} 
\cnode*(-3,11){.07}{a48} 
\cnode*(-2,11){.07}{a49} 
\cnode*(-2.5,11.5){.07}{a50} 
\cnode*(-1.5,11.5){.07}{a51} 
\cnode*(-2,12){.07}{a52} 
\pscircle(-2,12){.2}
\cnode*(-2,12.6){.07}{a53}
\cnode*(-2,13.2){.07}{a54} 
\cnode*(-2,13.8){.07}{a55}
\cnode*(-2,14.4){.07}{a56} 
\ncline{a1}{a2}
\ncline{a2}{a3}
\ncline{a3}{a4}
\ncline[linewidth=1.5pt]{a4}{a5}
\Bput{d}
\ncline{a5}{a6}
\ncline{a5}{a7}
\ncline{a6}{a8}
\ncline{a7}{a8}
\ncline{a7}{a9}
\ncline[linewidth=1.5pt]{a8}{a10}
\Bput{d}
\ncline{a8}{a11}
\ncline{a9}{a11}
\ncline{a10}{a12}
\ncline{a10}{a13}
\ncline[linewidth=1.5pt]{a11}{a13}
\Bput{d}
\ncline{a12}{a14}
\ncline{a12}{a15}
\ncline{a13}{a15}
\ncline{a13}{a16}
\ncline{a14}{a17}
\ncline{a14}{a18}
\ncline{a15}{a18}
\ncline{a15}{a19}
\ncline{a16}{a19}
\ncline{a17}{a20}
\ncline{a18}{a20}
\ncline{a18}{a21}
\ncline{a19}{a21}
\ncline[linewidth=1.5pt]{a19}{a22}
\Aput{d}
\ncline{a20}{a23}
\ncline{a21}{a23}
\ncline[linewidth=1.5pt]{a21}{a24}
\Aput{d}
\ncline{a22}{a24}
\ncline{a22}{a25}
\ncline[linewidth=1.5pt]{a23}{a26}
\Aput{d}
\ncline{a24}{a26}
\ncline{a24}{a27}
\ncline{a24}{a28}
\ncline{a25}{a28}
\ncline{a26}{a29}
\ncline{a26}{a30}
\ncline{a27}{a29}
\ncline{a27}{a31}
\ncline{a28}{a30}
\ncline{a28}{a31}
\ncline{a29}{a32}
\ncline{a29}{a33}
\ncline{a30}{a33}
\ncline{a31}{a33}
\ncline[linewidth=1.5pt]{a31}{a34}
\Aput{d}
\ncline{a32}{a35}
\ncline{a33}{a35}
\ncline[linewidth=1.5pt]{a33}{a36}
\Aput{d}
\ncline{a34}{a36}
\ncline{a34}{a37}
\ncline[linewidth=1.5pt]{a35}{a38}
\Aput{d}
\ncline{a36}{a38}
\ncline{a36}{a39}
\ncline{a37}{a39}
\ncline{a37}{a40}
\ncline{a38}{a41}
\ncline{a38}{a42}
\ncline{a39}{a42}
\ncline{a39}{a43}
\ncline{a40}{a43}
\ncline{a41}{a44}
\ncline{a42}{a44}
\ncline{a42}{a45}
\ncline{a43}{a45}
\ncline[linewidth=1.5pt]{a44}{a46}
\Bput{d}
\ncline{a44}{a47}
\ncline{a45}{a47}
\ncline{a46}{a48}
\ncline{a46}{a49}
\ncline[linewidth=1.5pt]{a47}{a49}
\Bput{d}
\ncline{a48}{a50}
\ncline{a49}{a50}
\ncline{a49}{a51}
\ncline{a50}{a52}
\ncline{a51}{a52}
\ncline[linewidth=1.5pt]{a52}{a53}
\Bput{d}
\ncline{a53}{a54}
\ncline{a54}{a55}
\ncline{a55}{a56}
\psset{labelsep=7pt}
\cnode(-1,14){.07}{c1} 
\uput[d](-1,14){-\g_1}
\cnode*(-.25,14){.07}{c2} 
\cnode*(.5,14){.07}{c3} 
\cnode*(1.25,14){.07}{c4} 
\uput[d](1.25,14){d}
\cnode*(1.25,14.75){.07}{c5} 
\cnode*(2,14){.07}{c6} 
\cnode*(2.75,14){.07}{c7} 
\cnode*(3.5,14){.07}{c8} 
\ncline[linestyle=dashed,dash=3pt 2pt]{c1}{c2}
\ncline{c2}{c3}
\ncline{c3}{c4}
\ncline{c5}{c4}
\ncline{c4}{c6}
\ncline{c6}{c7}
\ncline{c7}{c8}
\psline(3.38,13.88)(3.62,14.12)
\psline(3.38,14.12)(3.62,13.88)
\cnode*(2,1.9){.07}{b4} 
\pscircle(2,1.9){.2}
\cnode*(2,3.5){.07}{b8} 
\pscircle(2,3.5){.2}
\cnode*(2,4){.07}{b11} 
\pscircle(2,4){.2}
\cnode*(2,5.5){.07}{b19} 
\pscircle(2,5.5){.2}
\cnode*(2,6){.07}{b21} 
\pscircle(2,6){.2}
\cnode*(1.5,6.5){.07}{b23} 
\pscircle(1.5,6.5){.2}
\cnode*(2.5,7.5){.07}{b31} 
\pscircle(2.5,7.5){.2}
\cnode*(2,8){.07}{b33} 
\cnode*(2,8.5){.07}{b35} 
\cnode*(2,10){.07}{b44} 
\cnode*(2,10.5){.07}{b47} 
\cnode*(2,12){.07}{b52} 
\pscircle(2,12){.2}
\ncline[linestyle=dashed,dash=4pt 3pt]{b4}{b8}
\ncline{b8}{b11}
\ncline[linestyle=dashed,dash=4pt 3pt]{b11}{b19}
\ncline{b19}{b21}
\ncline{b21}{b23}
\ncline[linestyle=dashed,dash=4pt 3pt]{b21}{b31}
\ncline[linestyle=dashed,dash=4pt 3pt]{b23}{b33}
\ncline{b31}{b33}
\ncline{b33}{b35}
\ncline[linestyle=dashed,dash=4pt 3pt]{b35}{b44}
\ncline{b44}{b47}
\ncline[linestyle=dashed,dash=4pt 3pt]{b47}{b52}
\rput(3,1.9){\varnothing}
\cnode*(3,3.5){.07}{d1}
\psline(2.88,3.38)(3.12,3.62)
\psline(2.88,3.62)(3.12,3.38)
\cnode*(3,4){.07}{d2}
\cnode*(3.75,4){.07}{d3}
\ncline{d2}{d3}
\psline(3.63,3.87)(3.87,4.12)
\psline(3.63,4.12)(3.87,3.87)
\cnode*(3,5.5){.07}{d4}
\cnode*(3.75,5.5){.07}{d5}
\cnode*(4.5,5.5){.07}{d6}
\ncline{d4}{d5}
\ncline{d5}{d6}
\psline(4.38,5.38)(4.62,5.62)
\psline(4.38,5.62)(4.62,5.38)
\cnode*(3,6){.07}{d7}
\cnode*(3.75,6){.07}{d8}
\cnode*(4.5,6){.07}{d9}
\cnode*(5.25,6){.07}{d10}
\ncline{d7}{d8}
\ncline{d8}{d9}
\ncline{d9}{d10}
\psline(5.13,5.88)(5.37,6.12)
\psline(5.13,6.12)(5.37,5.88)
\cnode*(3,6.5){.07}{d11}
\cnode*(3.75,6.5){.07}{d12}
\cnode*(4.5,6.5){.07}{d13}
\cnode*(5.25,6.5){.07}{d14}
\cnode*(6,6.5){.07}{d15}
\ncline{d11}{d12}
\ncline{d12}{d13}
\ncline{d13}{d14}
\ncline{d14}{d15}
\psline(5.88,6.38)(6.12,6.62)
\psline(5.88,6.62)(6.12,6.38)
\cnode*(3.75,7.87){.07}{d16}
\cnode*(3.75,7.12){.07}{d17}
\cnode*(4.5,7.12){.07}{d18}
\cnode*(5.25,7.12){.07}{d19}
\cnode*(6,7.12){.07}{d20}
\ncline{d16}{d17}
\ncline{d17}{d18}
\ncline{d18}{d19}
\ncline{d19}{d20}
\psline(5.88,7)(6.12,7.24)
\psline(5.88,7.24)(6.12,7)
\cnode*(3,11.62){.07}{d21}
\cnode*(3.75,12.37){.07}{d22}
\cnode*(3.75,11.62){.07}{d23}
\cnode*(4.5,11.62){.07}{d24}
\cnode*(5.25,11.62){.07}{d25}
\cnode*(6,11.62){.07}{d26}
\ncline{d21}{d23}
\ncline{d22}{d23}
\ncline{d23}{d24}
\ncline{d24}{d25}
\ncline{d25}{d26}
\psline(5.88,11.5)(6.12,11.74)
\psline(5.88,11.74)(6.12,11.5)
$
\end{pspicture}
\caption{Kostant modules in a semi-regular block for $(E_7,E_6)$}
\label{F:E7}
\end{figure}

\begin{ex}
Let $\D=(E_7,E_6)$ and let $\O_S^\mu$ be a semi-regular block; i.e., the rank of the singular root system is $t=1$. The theorem says that the Kostant modules in $\O_S^\mu$ are parameterized by the subdiagrams of $\D^{(1)}=(D_6,D_5)$, since in the extended Dynkin diagram of $E_7$, $-\g_1$ is attached to the first simple root (see Figure~\ref{F:E7}). To describe explicitly the correspondence of the theorem, consider for example the case $J=\{d\}$. In the Hasse diagram of $\SW$ (left side of Figure~\ref{F:E7}) a node corresponds to an element of $\SWJ$ if it has an edge with label $d$ going up (these edges are indicated in bold). The right side of Figure~\ref{F:E7} shows the Hasse diagram of $\SWJ$, with dashed edges indicating coverings of length difference greater than one. This poset is isomorphic to the parabolic poset ${}^{S'}W'$ corresponding to $\D^{(1)}$. The subdiagrams of $\D^{(1)}$ corresponding to the Kostant modules are listed along the right side of the figure.
\end{ex}

\subsection{\it Minimal free resolutions of determinantal ideals\ }
Let $(G_{\RR},M_{\RR})$ be an irreducible Hermitian symmetric pair of 
non-compact type of rank $r$; i.e., $G_{\RR}$ is a simple non-compact (real) 
Lie group, $M_{\RR}$ is a maximal compact subgroup of $G_{\RR}$, and 
$G_{\RR}/M_{\RR}$ is a Hermitian symmetric space of rank $r$.
Let $\fg$ and $\fm$ be the complexified Lie algebras of $G_{\RR}$ and $M_{\RR}$
and let $\fg = \fu^{-}\oplus \fm \oplus \fu$ be the usual decomposition
of $\fg$ as an $\fm$-module. The complexification of the compact group 
$M_{\RR}$ acts on $\fu$ with finitely many orbits
and the closures of the  orbits form a chain of algebraic varieties
$X_{0}=\{0\} \subset X_{1}\subset \cdots \subset X_{r}=\fu$.
For the classical Hermitian symmetric pairs these varieties are determinantal
varieties. Let $I(X_k)$ denote the ideal of $X_k$ in $\CC[\fu]=S(\fu^{-})$;
i.e.,
$$
  I(X_{k})\ = \{ f\in \CC[\fu] \mid f(x) = 0 \quad \forall x \in X_k\}.
$$
It is well-known (see \cite{Jo:92}) that $I(X_k)$ is the annihilator of the
$k$th Wallach representation $L_{-kc\zeta}$ 
when looked upon as an $S(\fu^-)$-module. Here $L_{-kc\zeta}$ denotes the 
simple $\fg$-module with highest weight $-kc\zeta$, where $\zeta$ is the fundamental weight of $\fg$ that is perpendicular to all compact roots 
(i.e., perpendicular to $\Phi_{S}$) and 
$c$ is a constant that depends on the pair $(\fg,\fm)$.
Thus we can identify the coordinate ring $\CC[X_k]=\CC[\fu]/I(X_{k})$ with
$L_{-kc\zeta}$ as a graded $S(\fu^-)$-module (by sending $1\in \CC[X_k]$
to a highest weight vector in $L_{-kc\zeta}$). 
Moreover, the group $M_{\RR}$ acts naturally on $\CC[X_{k}]$ (since $M_{\RR}$
acts on $X_{k}$) and via this action, $\CC[X_k]\otimes F_{-kc\zeta}\simeq L_{-kc\zeta}$ as an $M_{\RR}$-module. (Note that 
$F_{-kc\zeta}$ is one-dimensional.)

Now suppose for simplicity that $\D$ is of simply laced type.
Then the module $L_{-kc\zeta}$ is singular for $0<k\leq r$, and $|J|=k$.
Furthermore, under the equivalence of categories $\O(\Phi,\Phi_S,\Phi_J) \simeq \O(\Phi',\Phi'_{S'},\varnothing)$, the Wallach representation $L_{-kc\zeta}$ 
corresponds to the finite dimensional module in the regular block
$\O(\Phi',\Phi'_{S'},\varnothing)$. Thus the BGG resolution
of $L_{-kc\zeta}$ is described by the poset ${}^{S'}W'$ corresponding
to $\D'=\D^{(k)}$. It was proved in \cite{EnHu:04a} that the BGG resolution
of $L_{-kc\zeta}$ gives a minimal free (and $M_{\RR}$-equivariant) resolution 
of $\CC[X_k]$ as an $R=\CC[\fu]$-module.

\begin{ex}
Consider the Hermitian symmetric pair corresponding to the diagram
$\D=(E_{7},E_{6})$. In this case, $M_{\CC}=E_{6}(\CC)$ and $\fu\simeq \CC^{27}$.
We would like to give the minimal free resolution of $\CC[X_{1}]$,
the coordinate ring of the closure of the minimal orbit. The corresponding 
Wallach representation is semi-regular and one finds that $J=\{d\}$.
Thus the minimal free resolution
of $\CC[X_{1}]$ is described by the poset ${}^{S}W^{J}\simeq {}^{S'}W'$
shown in  Fig.~\ref{F:E7}.  We will now give this resolution 
explicitly. Let $R=\CC[\fu]=S(\fu^{-})$ which we view as a graded ring.
For $m\in \ZZ$ define a graded free $R$-module, $R(m)$, by $R(m)_{i}=R_{i+m}$. 
Then the minimal free resolution of $\CC[X_{1}]$ is
\begin{equation*}
\begin{split}
0 \rightarrow & R(-15)\rightarrow R(-13)^{27}\rightarrow R(-12)^{78}
\rightarrow R(-10)^{351} \rightarrow R(-9)^{650}\rightarrow 
R(-8)^{351}\oplus R(-7)^{351}
\\
& \rightarrow R(-6)^{650} \rightarrow R(-5)^{351}
\rightarrow R(-3)^{78} \rightarrow R(-2)^{27}\rightarrow R \rightarrow 
\CC[X_{1}]\rightarrow 0.
\end{split}
\end{equation*}
Here the exponent of the term corresponding to 
$x\in  {}^{S}W^{J}$ is equal to the dimension of the $\fm$-module  $F_{x}$. 
We calculated these dimensions 
by using the Weyl dimension formula. 
The degree shifts appearing in the resolution can be directly obtained 
from Fig.~\ref{F:E7} as follows. The term corresponding to 
$x\in  {}^{S}W^{J}$ has a degree shift equal to the number of edges 
in a chain from $x$ to $w$ (in the poset ${}^{S}W^{J}$) {\it excluding\/} 
the edges corresponding to the singular root $d$. (Thus the degree shift 
is in general smaller than the difference $l(w)-l(x)$.)
\end{ex}

\section{Singular Blocks in the General Case} 
\Label{S:General}

\subsection{} 
Recall from our discussion in Section~\ref{S:HSsing} that for Hermitian symmetric pairs, every singular block is equivalent to a regular block (or a direct sum of two regular blocks) for some other Hermitian symmetric pair via equivalences of categories due to Enright and Shelton. 
However, in the non-Hermitian symmetric setting, singular blocks are 
more complicated than regular blocks. In this section we briefly explain how,
at least in principle, one can decide whether a simple module in a
singular block is a Kostant module via Kazhdan-Lusztig polynomials.
We then introduce a new ordering on the set $\SWJ$ that is equal to the Bruhat ordering when $J=\varnothing$ or $S=\varnothing$, but is different from the Bruhat ordering in general.

\subsection{\it KLV polynomials\ }

For regular blocks, it is possible to compute
$\Ext$-groups, and hence $\fu$-cohomology, using Kazhdan-Lusztig polynomials. 
It turns out that this is also possible for singular blocks.
Fix a singular block $\O_{S}^{\,\mu}$ and, as before, let $\SWJ$ denote the
set that parameterizes the simple modules in $\O_{S}^{\,\mu}$.
Then for $x,w\in \SWJ$ define the {\it Kazhdan-Lusztig-Vogan polynomial\/}  
${}^{S}\!P^{J}_{x,w}$ by the same formula as in Section \ref{SS:KLV}:
$$
 {}^{S}\!P^{J}_{x,w}=
 \sum_{i\geq 0} q^{\frac{l(w)-l(x)-i}{2}} \dim \Ext^{i}_{\O_{S}}(N_{x},L_{w}).
$$
There is a beautiful formula due to Soergel \cite{Soe:89} and  
Irving \cite{Irv:90a} that allows one to compute the ${}^{S}\!P^{J}_{x,w}$
in terms of (regular parabolic) KLV polynomials.
It says that for $x, w\in \SWJ$,
$$
  {}^{S}\!P^{J}_{x,w}=
  \sum_{z\in W_{J}} (-1)^{l(z)}\ {}^{S}\!P_{xz,w}\ .
$$
Note that if $x \in \SWJ$ and $z\in W_{J}$ then $l(xz)=l(x)+l(z)$. 
Hence the formula above is equivalent to 
$$
  \dim \Ext^{i}_{\O_{S}}(N_{w_{S}x\cdot\mu},L_{w_{S}w\cdot\mu})
  =\sum_{z\in W_{J}} (-1)^{l(z)} 
  \dim \Ext^{i-l(z)}_{\O_{S}}(N_{w_{S}xz\cdot\lambda},L_{w_{S}w\cdot\lambda}) 
$$
for all $i$, where $\lambda$ is any weight such that $\lambda+\rho$ is regular
and anti-dominant integral; e.g., $\lambda=-2\rho$. 
In \cite{Soe:89}, this formula is proved for singular blocks in 
ordinary category $\O$; i.e., for $S=\varnothing$. The general formula follows
from the following isomorphism (which is proved using a Lyndon-Hochschild-Serre 
spectral sequence argument as in \cite[Chap.\ 15]{EnSh:87})
$$
  \Ext^{i}_{\O_{S}}(N_{w_{S}y\cdot \mu}, L_{w_{S}w\cdot \mu}) \simeq
  \Ext^{i}_{\O}(M_{w_{S}y\cdot \mu}, L_{w_{S}w\cdot \mu}),
$$
where $M_{w_{S}y\cdot \mu}$ is the ordinary 
Verma module of highest weight $w_{S}y\cdot \mu$ in the block $\O^{\mu}$ 
of ordinary category $\O$.

\subsection{\it A new ordering on $\SWJ$.} \Label{SS:Extorder}
Using the formula for KLV polynomials ${}^{S}\!P^{J}_{x,w}$ from above, 
we started to identify Kostant modules in singular blocks for 
simple Lie algebras of small rank. 
We then noticed that in certain singular blocks, there exist 
simple modules $L_{w}$ that ``ought'' to be Kostant modules 
since all non-zero KLV polynomials ${}^{S}\!P^{J}_{x,w}=1$, but for
which the set $\{x\in \SWJ \mid {}^{S}\!P^{J}_{x,w}=1\}$ is not an interval of $\SWJ$. 
This led us to consider a different ordering on $\SWJ$. We propose to use this 
new ordering to obtain the rank function used to define Kostant modules 
(Definition \ref{D:Kostant}) for singular blocks of category $\O_S$.

\begin{defn}
For $x<w$ in $\SWJ$ define $\mu_{S}(x,w)$ by
$$\mu_{S}(x,w)=\dim \Ext^{1}_{\O_{S}} (L_{x},L_{w})
=\dim \Ext^{1}_{\O_{S}} (N_{x},L_{w}).
$$
For  $x, w\in \SWJ$, write $x\rightarrow_{\mu} w$ 
if $x<w$ in the Bruhat ordering, $\mu_{S}(x,w)\not=0$ 
and there is no $x<z<w$ in $\SWJ$ with $\mu_{S}(z,w)\not=0$.
Then define $\leq_{\mu}$ as the ordering on $\SWJ$ generated 
by the covering relations $x\rightarrow_{\mu} w$. We call this
ordering on $\SWJ$ the {\em $\mu$-ordering} or the {\em $\Ext^{1}$-ordering}.
\end{defn}

\begin{rem}
Note that by definition, $x\leq_{\mu} w$ implies $x\leq w$ in
the Bruhat ordering. If $J=\varnothing$ or $S=\varnothing$ then the 
$\Ext^{1}$-ordering and the Bruhat ordering coincide. The proof is easy
in the case when $J=\varnothing$. The case when $S=\varnothing$ follows
via Koszul duality for parabolic and singular category $\O$
(cf.\ \cite{Bac:99} and \cite{BGS:96}).
\end{rem}

We now give two examples to show the usefulness of the new ordering.

\subsection{\it Splitting of singular blocks for  Hermitian symmetric cases} \Label{SS:HSsplit}
In Section~\ref{S:HSsing} (cf.\ Table~\ref{T:HSsingdata}), we pointed out that that for $(\Phi,\Phi_{S})=(B_{n}, B_{n-1} )$ or $(C_{n}, A_{n-1})$, certain singular blocks $\O_{S}^{\, \mu}$ split further  into a direct sum 
$\O_{S}^{\, \mu}=\O_{1}\oplus \O_{2}$ with each of the summands equivalent 
to a regular block of another Hermitian symmetric pair. 
The $\Ext^{1}$-ordering reflects this nicely in the sense that the poset 
$\SWJ$ splits up into two disjoint
posets, each of which looks like a regular parabolic poset. 
On the other hand, the Bruhat ordering makes $\SWJ$ a connected poset. A simple module $L_w \in \O_{S}^{\, \mu}$ which correspond to a Kostant module in $\O_1$ or $\O_2$ ``ought'' to be a Kostant module, but in general the subset of $\SWJ$ parameterizing the cohomology of $L_w$ will only be an interval in the $\Ext^1$-ordering, not in the Bruhat ordering.

\subsection{\it A semi-regular block for $(F_4,C_{3})$}
\begin{figure}[ht]
\centering
\begin{pspicture}(-.3,-.35)(1,4.5)
\psset{linewidth=.5pt,labelsep=8pt,nodesep=0pt}
\small
$
\cnode*(-2,0){.07}{a1}
\pscircle(-2,0){.2}
\uput[r](-2,0){1}
\cnode*(-2,.75){.07}{a2}
\pscircle(-2,.75){.2}
\uput[r](-2,.75){2}
\cnode*(-2,1.5){.07}{a3} 
\pscircle(-2,1.5){.2}
\uput[r](-2,1.5){3}
\cnode*(-2,2.5){.07}{a4} 
\pscircle[linestyle=dashed,dash=4pt 2pt](-2,2.5){.2}
\uput[l](-2,2.5){?}
\uput[r](-2,2.5){4}
\cnode*(-2,3.25){.07}{a5} 
\uput[r](-2,3.25){5}
\cnode*(-2,4){.07}{a6} 
\pscircle[linestyle=dashed,dash=4pt 2pt](-2,4){.2}
\uput[l](-2,4){?}
\uput[r](-2,4){6}
\ncline{a1}{a2}
\ncline{a2}{a3}
\ncline[linestyle=dashed,dash=4pt 3pt]{a3}{a4}
\ncline{a4}{a5}
\ncline{a5}{a6}
\cnode*(2,0){.07}{b1} 
\pscircle(2,0){.2}
\uput[r](2,0){1}
\cnode*(2,.75){.07}{b2} 
\pscircle(2,.75){.2}
\uput[r](2,.75){2}
\cnode*(1.5,1.5){.07}{b3} 
\pscircle(1.5,1.5){.2}
\uput[l](1.5,1.5){3}
\cnode*(2.5,1.75){.07}{b4} 
\pscircle(2.5,1.75){.2}
\uput[r](2.5,1.75){4}
\cnode*(2,2.5){.07}{b5} 
\uput[r](2,2.5){5}
\cnode*(2,3.25){.07}{b6} 
\uput[r](2,3.25){6}
\pscircle(2,3.25){.2}
\ncline{b1}{b2}
\ncline{b2}{b3}
\ncline[linestyle=dashed,dash=4pt 3pt]{b2}{b4}
\ncline[linestyle=dashed,dash=4pt 3pt]{b3}{b5}
\ncline{b4}{b5}
\ncline{b5}{b6}
$
\end{pspicture}
\caption{The Bruhat and $\Ext^1$ posets for $(F_4, \{b,c,d\}, \{d\})$}
\label{F:F4}
\end{figure}
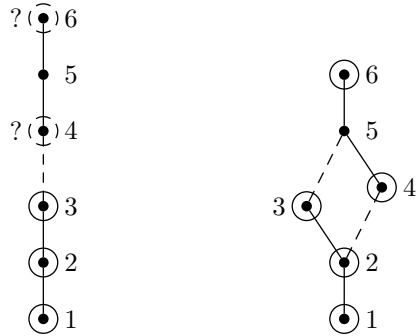

An interesting non-Hermitian symmetric example is obtained for 
$(\Phi,\Phi_{S})= (F_4,C_{3})$. Consider the category 
$\O(F_4, \{b,c,d\}, \{d\})$, where the labeling of the simple roots is 
the same as in Figure~\ref{F:F4C3}. The Hasse diagram of the poset 
$\SWJ$ with respect to the usual Bruhat ordering is given on the left side
of Figure~\ref{F:F4}. (Recall that in the Hasse diagram of ${}^{S}W$, Figure~\ref{F:F4C3},
a node that corresponds to an element of $\SWJ$ is a node 
having an edge with label $d$ going up from it.)
A calculation shows that $\Ext^1(V_3,L_4)=0$ and hence the 
$\Ext^{1}$-ordering on $\SWJ$ is different from the Bruhat ordering.
The Hasse diagram corresponding to the $\Ext^1$-poset is shown on the 
right side of Figure \ref{F:F4}. 
The structure of the category $\O(F_4, \{b,c,d\}, \{d\})$ is 
exactly the same as $\O(A_3,A_{1}\times A_{1},\varnothing)$. Not only are the $\Ext^1$ posets isomorphic, but also the radical filtrations of the generalized Verma modules in the two categories correspond:
$$
N_{1}= L_{1}
,\ 
N_{2}=
\mbox{\begin{tabular}{l}
$L_{2}$ \\
\hline
$L_{1}$\end{tabular}}\ ,\ 
N_{3}=\mbox{\begin{tabular}{c}
$L_{3}$ \\
\hline
$L_{2}$ \\
\end{tabular}}\ ,\ 
N_{4}=\mbox{\begin{tabular}{c}
$L_{4}$ \\
\hline
$L_{2}$ \\
\end{tabular}}\ ,\ 
N_{5}=\mbox{\begin{tabular}{c}
$L_{5}$ \\
\hline
$L_{4}$ $L_{3}$   $L_{1}$\\
\hline
$L_{2}$\\
\end{tabular}}\ ,\ 
N_{6}=\mbox{\begin{tabular}{c}
$L_{6}$ \\
\hline
$L_{5}$ \\
\hline
$L_{1}$\
\end{tabular}}\ .
$$
So, $H^{0}(\fu, L_{4})= F_{4}$, $H^{1}(\fu, L_{4})=F_{2}$, 
$H^{2}(\fu, L_{4})=F_{1}$, and $H^{i}(\fu, L_{4})=0$ for $i>2$.
Likewise, $H^{0}(\fu, L_{6})= F_{6}$, $H^{1}(\fu, L_{6})=F_{5}$, 
$H^{2}(\fu, L_{6})=F_{3}\oplus F_{4}$, $H^{3}(\fu, L_{6})=F_{2}$, 
$H^{4}(\fu, L_{6})=F_{1}$, and $H^{i}(\fu, L_{4})=0$ for $i>4$.
In particular, using the rank function of the $\Ext^1$-poset, $L_4$  and $L_6$ become Kostant modules.

\let\section=\oldsection

\begin{thebibliography}{Deo2}

\bibitem[Bac]{Bac:99}
Erik Backelin, {\em Koszul duality for parabolic and singular category {$\scr
  O$}}, Represent. Theory \textbf{3} (1999), 139--152 (electronic).
  
\bibitem[BGG]{BGG:75}  
Joseph N. Bernstein, Israel M., Gelfand, and Sergei~I. Gelfand,{\em Differential operators on the base affine space and a study of ${\germ g}$-modules}, Lie groups and their representations (Proc. Summer School, Bolyai J\'anos Math. Soc., Budapest, 1971), pp. 21--64. Halsted, New York, 1975.

\bibitem[BGS]{BGS:96}
Alexander Beilinson, Victor Ginzburg, and Wolfgang Soergel, {\em Koszul duality
  patterns in representation theory}, J. Amer. Math. Soc. \textbf{9} (1996),
  no.~2, 473--527.

\bibitem[BL]{BiLa:00}
Sara Billey and V.~Lakshmibai, {\em Singular loci of {S}chubert varieties},
  Birkh\"auser, Boston, 2000.

\bibitem[BiP]{BiPo:05}
Sara Billey and Alexander Postnikov, {\em Smoothness of {S}chubert varieties
  via patterns in root subsystems}, Adv. in Appl. Math. \textbf{34} (2005),
  no.~3, 447--466.

\bibitem[BN]{BoNa:05}
Brian~D. Boe and Daniel~K. Nakano, {\em Representation type of the blocks of
  category {${\scr O}\sb S$}}, Adv. Math. \textbf{196} (2005), no.~1, 193--256.

\bibitem[Bou]{Bou:68}
N.~Bourbaki, {\em \'{E}l\'ements de math\'ematique. {F}asc. {XXXIV}. {G}roupes
  et alg\`ebres de {L}ie. {C}hapitre {IV}: {G}roupes de {C}oxeter et syst\`emes
  de {T}its. {C}hapitre {V}: {G}roupes engendr\'es par des r\'eflexions.
  {C}hapitre {VI}: syst\`emes de racines}, Actualit\'es Scientifiques et
  Industrielles, No. 1337, Hermann, Paris, 1968.

\bibitem[BrP]{BrPo:99}
Michel Brion and Patrick Polo, {\em Generic singularities of certain {S}chubert
  varieties}, Math. Z. \textbf{231} (1999), no.~2, 301--324.

\bibitem[Car]{Car:94}
James~B. Carrell, {\em The {B}ruhat graph of a {C}oxeter group, a conjecture of
  {D}eodhar, and rational smoothness of {S}chubert varieties}, Proc.\ Sympos.\
  Pure Math., vol. 56, part 1, pp.~53--61, Amer. Math. Soc., Providence, RI,
  1994.

\bibitem[CK]{CaKu:03}
James~B. Carrell and Jochen Kuttler, {\em Smooth points of {$T$}-stable
  varieties in {$G/B$} and the {P}eterson map}, Invent. Math. \textbf{151}
  (2003), no.~2, 353--379.

\bibitem[CC]{CaCo:87}
Luis~G. Casian and David~H. Collingwood, {\em The {K}azhdan-{L}usztig
  conjecture for generalized {V}erma modules}, Math. Z. \textbf{195} (1987),
  no.~4, 581--600.

\bibitem[Col]{Col:85}
David~H. Collingwood, {\em The {${\germ n}$}-homology of {H}arish-{C}handra
  modules: generalizing a theorem of {K}ostant}, Math. Ann. \textbf{272}
  (1985), no.~2, 161--187.

\bibitem[Deo1]{Deo:77}
Vinay~V. Deodhar, {\em Some characterizations of {B}ruhat ordering on a
  {C}oxeter group and determination of the relative {M}{\"o}bius function},
  Invent. Math. \textbf{39} (1977), 187--198.

\bibitem[Deo2]{Deo:87}
Vinay~V. Deodhar, {\em On some geometric aspects of {B}ruhat orderings. {II}.
  {T}he parabolic analogue of {K}azhdan-{L}usztig polynomials}, J. Algebra
  \textbf{111} (1987), no.~2, 483--506.

\bibitem[Enr]{Enr:88}
Thomas~J. Enright, {\em Analogues of {K}ostant's {${\germ u}$}-cohomology
  formulas for unitary highest weight modules}, J. Reine Angew. Math.
  \textbf{392} (1988), 27--36.

\bibitem[EH1]{EnHu:04a}
Thomas~J. Enright and Markus Hunziker, {\em Resolutions and Hilbert series 
of determinantal varieties and unitary highest weight modules},
J.~Algebra \textbf{273} (2004),  608--639.

\bibitem[EH2]{EnHu:04b}
\bysame {\em Resolutions and Hilbert series of the unitary highest weight modules 
of the exceptional groups},  Represent. Theory \textbf{8} (2004), 15--51.

\bibitem[ES1]{EnSh:87}
Thomas~J. Enright and Brad Shelton, {\em Categories of highest weight modules:
  applications to classical {H}ermitian symmetric pairs}, Mem. Amer. Math. Soc.
  \textbf{67} (1987), no.~367, iv+94.

\bibitem[ES2]{EnSh:89}
\bysame, {\em Highest weight modules for {H}ermitian symmetric pairs of
  exceptional type}, Proc. Amer. Math. Soc. \textbf{106} (1989), no.~3,
  807--819.

\bibitem[Hum]{Hum:90}
James~E. Humphreys, {\em Reflection groups and {C}oxeter groups}, Cambridge
  studies in advanced mathematics, vol.~29, Cambridge University Press,
  Cambridge, 1990.

\bibitem[Irv]{Irv:90a}
Ronald~S. Irving, {\em Singular blocks of the category {${\mathcal O}$}}, Math.
  Z. \textbf{204} (1990), 209--224.
  
\bibitem[Jo]{Jo:92}  
Anthony Joseph, {\em Annihilators and associated varieties of unitary 
highest weight modules},
Ann. scient. Ec. Norm. Sup., $4^e$ s\'erie, \textbf{25} (1992), 1--45.


\bibitem[KL]{KaLu:79}
David Kazhdan and George Lusztig, {\em Representations of {C}oxeter groups and
  {H}ecke algebras}, Invent. Math. \textbf{53} (1979), 165--184.

\bibitem[Kos]{Kos:61}
Bertram Kostant, {\em Lie algebra cohomology and the generalized {B}orel-{W}eil
  theorem}, Ann. of Math. \textbf{74} (1961), 329--387.

\bibitem[Lep]{Lep:77}
J.~Lepowsky, {\em A generalization of the {B}ernstein-{G}elfand-{G}elfand
  resolution}, J. Algebra \textbf{49} (1977), no.~2, 496--511.

\bibitem[Pro]{Proc:82}
Robert~A. Proctor, {\em Classical {B}ruhat orders and lexicographic
  shellability}, J. Algebra \textbf{77} (1982), 104--126.

\bibitem[RC]{Roc:80}
Alvany Rocha-Caridi, {\em Splitting criteria for {${\germ g}$}-modules induced
  from a parabolic and the {B}er\v nste\u\i n-{G}el\cprime fand-{G}el\cprime
  fand resolution of a finite-dimensional, irreducible {${\germ g}$}-module},
  Trans. Amer. Math. Soc. \textbf{262} (1980), no.~2, 335--366.

\bibitem[Sch]{Sch:81}
Wilfried Schmid, {\em Vanishing theorems for {L}ie algebra cohomology and the
  cohomology of discrete subgroups of semisimple {L}ie groups}, Adv. in Math.
  \textbf{41} (1981), no.~1, 78--113.

\bibitem[Soe]{Soe:89}
W.~Soergel, {\em {$\germ n$}-cohomology of simple highest weight modules on
  walls and purity}, Invent. Math. \textbf{98} (1989), no.~3, 565--580.

\end{thebibliography}

\def\germ{\mathfrak}\def\cprime{$'$}\def\scr{\mathcal}

\end{document}